\documentclass[a4paper, 11 pt , english]{article}

\usepackage{ulem}
\usepackage[english]{babel}
\usepackage{amsthm}
\usepackage{pstricks,pst-node,pst-tree}
\usepackage{amssymb}
\usepackage[utf8]{inputenc}
\usepackage[T1]{fontenc} 
\usepackage{fancybox} 
\usepackage{alltt}
\usepackage{graphicx}
\usepackage{bm}
\usepackage{lmodern}           
\usepackage{babel}
\usepackage[top = 3cm, bottom = 3 cm, right = 3cm, left = 3 cm]{geometry}
\usepackage{appendix}
\usepackage{float}
\usepackage{xcolor}
\usepackage{fancyhdr}
\usepackage{tikz}
\usetikzlibrary{arrows.meta}
\usepackage{csquotes}
\usepackage{amsmath, amssymb, mathtools, stmaryrd, physics, empheq, amsthm}

\usepackage{graphicx}
\graphicspath{{images/}}
\usepackage{caption}
\usepackage{subcaption}
\usepackage{float}
\usepackage{xcolor}
\usepackage{tikz}
\usetikzlibrary{arrows}
\usepackage{pgfplots}
\pgfplotsset{compat=1.15}
\usepackage[colorlinks,hyperindex,bookmarks,linkcolor=blue,citecolor=blue,urlcolor=blue]{hyperref}

\usepackage[linesnumbered, ruled]{algorithm2e}

\usepackage{enumitem}
\usepackage{dsfont}
\usepackage[most]{tcolorbox}
\usepackage{appendix}

\newtheorem{Theorem}{Theorem}[section]

\newtheorem{Proposition}[Theorem]{Proposition}
\newtheorem{Assumption}[Theorem]{Assumption}

\newtheorem{Lemma}[Theorem]{Lemma}
\newtheorem{Corollary}[Theorem]{Corollary}
\newtheorem{Remark}[Theorem]{Remark}

\def \be{\begin{eqnarray}}
\def \ee{\end{eqnarray}}
\def \b*{\begin{align*}}
\def \e*{\end{align*}}

\def \E{\mathbb{E}}
\def \F{\mathbb{F}}

\def \H{\mathbb{H}}
\def \L{\mathbb{L}}

\def \P{\mathbb{P}}

\def \R{\mathbb{R}}
\def \S{\mathbb{S}}

\def \V{\mathbb{V}}

\def \]{[\,\!\![}
\def \]{]\,\!\!]}

\def \proof{{\noindent \bf Proof. }}

\def \1{{\bf 1}}

\def \ep{\hbox{ }\hfill$\Box$}

\def\Ac{{\cal A}}
\def\Bc{{\cal B}}

\def\Ec{{\cal E}}
\def\Fc{{\cal F}}

\def\Hc{{\cal H}}

\def\Lc{{\cal L}}

\def\Kc{{\cal K}}

\def\Tc{{\cal T}}

\def\Xc{{\cal X}}

\def\namedlabel#1#2{\begingroup
  #2%
  \def\@currentlabel{#2}%
  \phantomsection\label{#1}\endgroup
}
\newcommand{\vertiii}[1]{{\left\vert\kern-0.25ex\left\vert\kern-0.25ex\left\vert #1 
 \right\vert\kern-0.25ex\right\vert\kern-0.25ex\right\vert}}

\def\eps{\varepsilon}

\newcommand{\tl}{\tilde}
\newcommand{\lm}{\lambda}
\newcommand{\gm}{\gamma}
\let\T\undefined
\newcommand{\T}{\Tc}
\newcommand{\dU}{\delta U}
\newcommand{\dV}{\delta V}

\let\L\undefined
\newcommand{\L}{\mathbb{L}}
\let\H\undefined
\newcommand{\H}{\mathbb{H}}
\newcommand{\Hb}{\mathbb{H}^2_\mathrm{BMO}}
\newcommand{\Hbtau}{\mathbb{H}^2_{\mathrm{BMO},\tau}}
\DeclareMathOperator*{\essinf}{ {\rm essinf}}
\DeclareMathOperator*{\esssup}{ {\rm esssup}}
\DeclareMathOperator*{\argmin}{ {\rm argmin}}
\DeclareMathOperator*{\argmax}{ {\rm argmax}}

\newlist{romanenumerate}{enumerate}{1}
\setlist[romanenumerate,1]{label=\textup{(\roman*)}, leftmargin=2em}

\usepackage[backend=biber,style=numeric, maxcitenames=99,   
 maxbibnames=99    
 ]{biblatex}

\usepackage{wrapfig}
\usepackage{epsf}

\DeclareFieldFormat{title}{#1}
\DeclareFieldFormat[article]{title}{#1}
\DeclareFieldFormat[inproceedings]{title}{#1}
\DeclareFieldFormat[book]{title}{#1}
\DeclareFieldFormat[incollection]{title}{#1}

\DeclareFieldFormat[article]{journaltitle}{\textit{#1}}
\DeclareFieldFormat[inproceedings]{booktitle}{\textit{#1}}
\DeclareFieldFormat[incollection]{booktitle}{\textit{#1}}
\DeclareFieldFormat[book]{title}{\textit{#1}}

\AtEveryBibitem{%
 \ifentrytype{article}{\clearfield{publisher}}{}
}

\AtEveryBibitem{%
 \iffieldundef{doi}{}{ \clearfield{url} }
}

\AtBeginBibliography{%
 \sloppy
 \emergencystretch=3em
}

\title{\bf Sensitivity Analysis of Distributionally Robust
BSDEs and RBSDEs}
\author
{     Arthur Compoint\thanks{Ecole Polytechnique, CMAP/
                      Tandon School of Engineering, 
                      arthur.compoint@polytechnique.edu. }
                      \and 
      Nathan Sauldubois\thanks{New York University, Tandon School of Engineering
            ns6982@nyu.edu}\and 
       Nizar Touzi\thanks{New York University, Tandon School of Engineering,
nizar.touzi@nyu.edu. Nizar Touzi gratefully acknowledges financial support
from the National Science Foundation under grant DMS-2508581.}
      }

\date{\today}  

\begin{document}

\maketitle

\vspace{-8mm}

\begin{abstract}	
We examine the sensitivity properties of backward stochastic differential equations and reflected backward stochastic differential equations, which naturally arise in the context of optimal control and optimal stopping problems. 
Motivated by sensitivity analysis questions in distributionally robust optimization (DRO) control and optimal stopping problems, we establish explicit formulas for the corresponding sensitivities under drift perturbations of a reference measure. 
Our work is closely related to \citeauthor{bartl2023sensitivity} \cite{bartl2023sensitivity}.
In contrast to the existing literature, our analysis is carried out within a general non-Markovian framework.
\end{abstract}


\vspace{1mm}

\vspace{-4mm}
\section{Introduction}

\vspace{-2mm}

In the field of stochastic analysis, backward stochastic differential equations (BSDEs), introduced by \citeauthor{pardoux1990adapted} \cite{pardoux1990adapted}, and their reflected counterparts, first studied by \citeauthor{el1997reflected} \cite{el1997reflected}, are fundamental tools for addressing optimal control and optimal stopping problems in the non-Markovian setting. 
While closely connected to the classical partial differential equation (PDE) framework, these equations offer greater flexibility by naturally accommodating non-Markovian structures. 
From the perspective of applications, BSDEs and RBSDEs are of particular relevance in mathematical finance, where optimal control plays a central role. 
Examples include hedging under market imperfections (see Example $1.2$ in \citeauthor{el1997backward} \cite{el1997backward}), portfolio optimization as studied in \citeauthor{merton1971optimum} \cite{merton1971optimum}, \citeauthor{rouge2000pricing} \cite{rouge2000pricing} and \citeauthor{hu2005utility} \cite{hu2005utility} in the non-Markovian case, the pricing of American options and optimal stopping problems, see \citeauthor{el1997reflected} \cite{el1997reflected} and Dynkin games as in \citeauthor{cvitanic1996backward} \cite{cvitanic1996backward}. 
However, these approaches rely on the existence of a model to describe the behavior of the quantity of interest. 

Our objective in this paper is to introduce the notion of model risk through distributionally robust optimization (DRO) in the context of backward SDEs and reflected backward SDEs.
This risk is traditionally assessed by sensitivity measures such as the Greeks (see the chapter on Greeks in \citeauthor{hull1997options} \cite{hull1997options}).
However, Greeks may underestimate risk, since they only account for deviations within a given family of models.
An alternative approach, which aims to mitigate risk, relies on models of uncertain volatility, as developed by \citeauthor{lyons1995uncertain} \cite{lyons1995uncertain} and \citeauthor{avellaneda1995pricing} \cite{avellaneda1995pricing}, where the key idea is to construct hedging strategies robust to worst-case scenarios.
The idea of mitigating risk by anticipating worst-case scenarios is commonly referred to as the Knightian approach, developed by \citeauthor{knight1921risk} \cite{knight1921risk}, and has been extensively studied in the distributionally robust optimization literature. 
In DRO, the model-misspecification is modeled through a deviation set $D$ of possible models. 
An early instance of this idea appears in \citeauthor{scarf1958min} \cite{scarf1958min}, where $D$ is defined through a set of models satisfying moment constraints.
More recently, the deviation set has typically been defined as a ball for a divergence or a distance on the space of probability measures (see \citeauthor{ben2013robust} \cite{ben2013robust}, \citeauthor{hu2013kullback} \cite{hu2013kullback}, for divergence-based criteria or \citeauthor{mohajerin_esfahani_data-driven_2018} \cite{mohajerin_esfahani_data-driven_2018}, \citeauthor{blanchet2019quantifying} \cite{blanchet2019quantifying} for optimal-transport-based criteria) 

In decision theory and macroeconomics, robust control (closely related to DRO) has been extensively studied as a framework to address model misspecification. In \citeauthor{hansen1999robust} \cite{hansen1999robust}, \citeauthor{hansen2001robust} \cite{hansen2001robust}, \citeauthor{anderson2003quartet} \cite{anderson2003quartet}, \citeauthor{hansen2006robust} \cite{hansen2006robust}, and \citeauthor{hansen2024risk} \cite{hansen2024risk}, the ambiguity set is constructed using a relative entropy criterion. The central idea of robust control, or more generally DRO, is that the decision maker interacts with an adversarial agent who seeks to minimize her gains by selecting the worst-case model within a prescribed set of probability measures (see \citeauthor{chen2002ambiguity} \cite{chen2002ambiguity}; \citeauthor{gilboa1989maxmin} \cite{gilboa1989maxmin}).

DRO applied to stochastic control problems has been studied in various contexts. For instance, \citeauthor{hansen2001robust} \cite{hansen2001robust} analyze robust control under an entropic ambiguity criterion. More recently, \citeauthor{bartl2023sensitivity} \cite{bartl2023sensitivity} investigated a Markovian optimization problem, addressing both drift and volatility uncertainty to derive first-order sensitivity results. Building on this line of work, \citeauthor{jiang2024sensitivity} \cite{jiang2024sensitivity} studied the sensitivity of path functionals under adapted Wasserstein deviations.

{\bf Our contribution} Inspired by the aforementioned works, we study DRO under $\L^\infty$ and $\L^2$ model perturbations of the drift, but in a non-Markovian framework. 
Let us illustrate the main results of the paper in the following simple setting. 
Let $A \subset \R^d$ be a bounded set and denote by $\Ac$ the class of $A$-valued admissible controls. 
We define $\E^\alpha := \E^{\P^\alpha}$, where $\P^\alpha$ is the probability measure under which the canonical process $X$ has drift $\alpha$.
We introduce the distributionally robust control problems,
\vspace{-3mm}
\begin{equation}\label{eqdef:sensi_opti_control}
 \overline{V}_\infty(r):= \inf_{\alpha\in\Ac}\sup_{\beta\in\Bc^\infty(r)}\E^{{\alpha+\beta}} [ \xi ] \, \, \text{and} \, \, 
 \overline{V}_2(r):= \inf_{\alpha\in\Ac}\sup_{\beta\in \Bc^2 (r)}\E^{{\alpha+\beta}} [\xi ],
 \vspace{-2mm}
\end{equation}
where $\Bc^\infty (r)$ and $\Bc^2(r)$ are $\L^\infty$ and $\L^2$ deviations. 
We prove that, under suitable conditions, $\overline{V}_\infty$ and $\overline{V}_2$ are differentiable at the origin and
\vspace{-2mm}
\begin{equation*}
\overline{V}_\infty'(0) \! =
\!
 \Vert Z \Vert_{\L^1 (\P^{\alpha^*})}:=
\E^{\P^{\alpha^*}}
[ {\textstyle \int}_0^T
\vert Z_t \vert \mathrm{d}t
]
\, \, \text{and} \, \, 
\overline{V}_2'(0) = \Vert Z \Vert_{\L^2 (\P^{\alpha^*})}:=\E^{\P^{\alpha^*}}[
{\textstyle \int}_0^T\vert Z_t \vert^2\mathrm{d}t
 ]^\frac12,
 \vspace{-3mm}
\end{equation*}
where $f_t(Y_t,Z_t):=\inf_{a \in A }\{a \cdot Z_t \}$, $\alpha^*_t$ 
is a corresponding minimizer,
and $(Y,Z)$ is the solution of the BSDE
\vspace{-3.5mm}
$$
 \mathrm{d} Y_t=Z_t\cdot \mathrm{d} X_t- f_t(Y_t,Z_t)\mathrm{d}t \,\, \text{and} \,\,
 Y_T=\xi.
 \vspace{-3mm}
$$

\noindent These results are closely related to those obtained by \citeauthor{bartl2023sensitivity} \cite{bartl2023sensitivity} and \citeauthor{jiang2024sensitivity} \cite{jiang2024sensitivity}. 
In particular, \citeauthor{jiang2024sensitivity} \cite{jiang2024sensitivity} established that the first-order sensitivity of a continuous-time model is given by the Malliavin derivative of the functional. 
Our findings align with this perspective, since the Malliavin derivative of $Y$ coincides with $Z$, as demonstrated in Proposition $5.3$ \citeauthor{el1997backward} \cite{el1997backward}.

\noindent We further extend the analysis to the setting of mixed optimal control/stopping problems, under $\L^\infty$ and $\L^2$ drift perturbations as illustrated by the problems:
\vspace{-2mm}
\begin{equation}\label{eqdef:sensi_optima_stopping}
\overline{\V}_\infty(r):=\inf_{\tau\in\Tc} \inf_{\alpha\in\Ac} \sup_{\beta\in\Bc^\infty(r)}\E^{\alpha + \beta}[\xi_\tau]
\,\, \text{and} \,\, 
\overline{\V}_2(r):= \inf_{\alpha\in\Ac} \inf_{\tau\in\Tc}\sup_{\beta\in \Bc^2 (r)}\E^{ \alpha + \beta}[\xi_\tau]
,\vspace{-2mm}
\end{equation}
where $(\xi_t)_{t\in[0,T]}\in\H^0$ is an obstacle process.
Under certain assumptions, we show that both $\overline{\V}_\infty$ and $\overline{\V}_2$ are differentiable at $0$ and compute their derivatives, which are given by
$
 \overline{\V}_\infty'(0)= \Vert Z\mathbf{1}_{[0,\tl{\tau}]} \Vert_{\L^1 (\P^{\alpha^*})}
$ and $
 \overline{\V}_2'(0)= \Vert Z\mathbf{1}_{[0, \tl{\tau}]} \Vert_{\L^2 (\P^{\alpha^*})},
$
where
$(Y,Z,K)$ is the solution of the RBSDE
\vspace{-5mm}
$$
\mathrm{d} Y_t = -f_t(Y_t, Z_t)\mathrm{d}t + Z_t\cdot \mathrm{d} X_t+\mathrm{d} K_t,
\quad Y_T=\xi_T
\quad
  Y\leq\xi, \quad\int_0^T(Y_t-\xi_t)\mathrm{d} K_t=0,
  \vspace{-2.5mm}
$$
$\tl{\tau}$ is the first hitting time $\tl{\tau}:=\inf\{t\geq0:Y_t=\xi_t\}$ and $f_t(Y_t,Z_t):=\inf_{a \in A }\{a \cdot Z_t \}$. 
The paper is organized as follows. We first derive an $\L^\infty$ sensitivity formula for BSDEs. We then study both $\L^2$ and $\L^\infty$ sensitivities for a general non-Markovian control problem, and subsequently for a mixed optimal stopping/control problem. We also compare our results with related works in the literature. Finally, we provide numerical illustrations in the case of a European lookback option under different borrowing and lending rates, and of an American lookback option problem.
\vspace{-2mm}

\section{Notations and Definitions}
Let $\Omega:=\mathcal{C}^0([0,T],\R^d)$, endowed with its canonical process $X=(X_t)_{t\in[0,T]}$, $X_t(\omega)=\omega(t)$. For a measurable set $A\subset\R^d$, let $\Bc_A$ denote its Borel $\sigma$-field. Let $\F=(\Fc_t)_{t\in[0,T]}$ be a filtration satisfying the usual conditions, and let $\P^0$ be the Wiener measure on $(\Omega,\F)$, under which $X$ is a Brownian motion. Throughout the paper, $\E=\E^0$ denotes expectation under $\P^0$ when no ambiguity arises, and $\E_t[\cdot]:=\E[\cdot\mid\Fc_t]$. We denote by $\H^0$ the set of progressively measurable processes.
Define
\vspace{-1.5mm}
\begin{itemize}
\item $\Tc_s^t$ is the set of $[s,t]$-valued $\F$-stopping times. We set
\vspace{-2.5mm}
  \begin{equation}\label{eqdef:Stopping time}
    \Tc:=\Tc_0^T \,\, ,\,\, \Tc_t:=\Tc_t^T.
    \vspace{-4.5mm}
  \end{equation}
  \item We write $\L^p := \L^p(\Fc_T):=\{Y\in\L^0(\Fc_T): \Vert Y \Vert_{\L^p}^p:=\E[ \vert Y \vert ^p]<\infty \}$ when no confusion is possible.
  \vspace{-3.5mm}
  \item For $k\in\L^\infty$, $\L^p(k):=\{Y\in\H^0: \Vert Y \Vert_{\L^p(k)}^p:=\E [\int_0^Te^{-\int_0^tk_s\mathrm{d}s} \vert Y_t \vert ^p\mathrm{d}t ]<\infty \}$.
  \vspace{-3mm}
  \item $\H^{q}
  \! := \! \big\{Y \! \in \! \H^0 \!:\! \Vert Y \Vert_{\H^{q}}^q \!:= \!\E\big[\big(\int_0^T \vert Y_t \vert ^2\mathrm{d}t\big)^{\frac{q}{2}}\big] \! < \! \infty \! \big\}$ and 
    \vspace{-3mm}
  \item $
\H^{q}_{\mathrm{loc}}\!
:=\!
\Big\{\!
Y \!\! \in \H^{0}
\!:\! 
\exists (\tau_n)_{n\ge1} \in \Tc_0^T,\,
\tau_n \! \uparrow \! T \ \text{a.s.},\ 
\text{and }
\mathbf{1}_{[0,\tau_n]} Y \in \H^{q}
\ \text{for all } n \ge 1
\Big\}.
$
  \vspace{-3mm}
  \item $\S^p:= \{Y\in\H^0:Y\text{ is continuous \textit{a.s.}, } \Vert Y \Vert_{\S^p}^p:=\E [\sup_{0\leq t\leq T} \vert Y_t \vert ^p ]<\infty \}$.
  \vspace{-3mm}
  \item $\mathbb{A}^p:= \{Y\in\S^p:Y\text{ is nondecreasing \textit{a.s.}, } Y_0=0 \}$.
  \vspace{-3mm}
  \item $\L^\infty$ denotes the space of essentially bounded random variables or stochastic processes, depending on the context. 
  \vspace{-2mm}
  \item $\S^\infty:=\big\{ Y\in \H^0: Y\text{ continuous, \textit{a.s.} and } \, \, \sup_{ 0 \leq t \leq T} \vert Y_t \vert \in \L^\infty \big\}$.
  \vspace{-2mm}
  \item $\H^p_\textrm{BMO} := \{Z\in\H^0: \Vert Z \Vert_{\H^p_\textrm{BMO}}^p:= \sup_{\tau\in\Tc} \Vert \E [\int_\tau^T \vert Z_t \vert ^p\mathrm{d}t \vert \Fc_\tau ]\Vert _\infty<\infty \}$.
By Chapter $2$ of \citeauthor{kazamaki2006continuous} \cite{kazamaki2006continuous}, $\Hb\subset\cap_{p\geq1}\H^p$. 
More precisely, for all $p\geq1$, there exists $C_p>0$ such that, for all $\beta\in\Hb$, 
\vspace{-2mm}
  \begin{equation}\label{ineq:HBMOLp} 
  \Vert \beta \Vert_{\H^p}\leq C_p \Vert \beta \Vert_{\Hb}.
  \end{equation}
  \end{itemize}
\vspace{-4mm}

\noindent When a random variable or a process takes values in a vector space, we indicate it explicitly, \textit{e.g.}, by writing $\H^p(\R^d)$. 
Furthermore, when needed, we specify the underlying probability measure, and we add an index $\tau$ whenever we are only interested in the process up to $\tau\in\Tc$. 
For instance, when $k\in\L^\infty$, we write
\vspace{-4mm}
\begin{equation*}
\L^p_\tau(\mathbb{Q},k):= \big\{Y\in\H^0: \Vert Y \Vert_{\L^p_\tau(\mathbb{Q},k)}^p:=
\E^\mathbb{Q} \Big[\int_0^{\tau}e^{-\int_0^tk_s\mathrm{d}s} \vert Y_t \vert ^p\mathrm{d}t \Big]<\infty \big\}.
\vspace{-3mm}
\end{equation*}
For $\beta\in\H^2_{\text{loc}}(\R^d)$, we define $\Ec(\int_s^t\beta_u\cdot \mathrm{d} X_u):= {\rm exp} (\int_s^t\beta_u\cdot \mathrm{d} X_u-\frac12\int_s^t \vert \beta_u \vert ^2\mathrm{d}u)$, the Doléans-Dade exponential of $\beta$. 
Whenever $\Ec(\int_0^T\beta_u\cdot \mathrm{d} X_u)$ is a martingale, we define the probability measure $\P^\beta$, equivalent to $\P^0$, by its density $\frac{d\P^\beta}{d\P^0}=\Ec(\int_0^T\beta_u\cdot \mathrm{d} X_u)$, so that $W^\beta:=X-\int_0^\cdot\beta_s\mathrm{d}s$ is a Brownian motion under $\P^\beta$. 
To simplify notation, we write $\E^\beta:=\E^{\P^\beta}$. 
Now, for a generator $f$, a progressively measurable process $(\xi_t)_{t \leq T}$, and a bounded stopping time $\tau$, we define the backward stochastic differential equation
\vspace{-1.5mm}
\begin{equation}\label{eqdef:BSDE random term time}
{\rm BSDE}_\tau(f, \xi)  
\,\, \,\, \mathrm{d} Y_t = Z_t\cdot \mathrm{d} X_t - f_t( Y_t, Z_t) \mathrm{d} t \,\, , \,\, Y_\tau =\xi_\tau
.
\vspace{-1.5mm}
\end{equation}
For simplicity, when $\tau = T$ and $\xi$ is an $\Fc_T$-measurable random variable, we write ${\rm BSDE}(f, \xi)$. 
Finally, for a generator $f$ and a barrier process $(\xi_t)_{0 \leq t \leq T}$, we define the reflected BSDE, which is a standard BSDE with the addition of a non-decreasing process $K$, by
\vspace{-1mm}
\begin{equation*}
{\rm RBSDE}(f, \xi) \,\, 
  \mathrm{d} Y_t = Z_t\cdot \mathrm{d} X_t+ \mathrm{d} K_t - f(Y_t, Z_t) \mathrm{d}t , \,\, Y_T=\xi_T, \,\, 
  Y\leq\xi, \,\, \int_0^T(Y_t-\xi_t)\mathrm{d} K_t=0
.\end{equation*}

\section{Main Results}

We first study the sensitivity of a general BSDE with respect to $\L^\infty$ deviations. This result will be used in the next subsections to derive sensitivity formulas for optimal control problems.
We then investigate the differentiability at $0$ of a class of problems that includes \eqref{eqdef:sensi_opti_control}. 
Finally, we address the sensitivity of mixed optimal control/stopping problems, thereby extending our analysis to a class of problems that includes \eqref{eqdef:sensi_optima_stopping}.

\subsection{Sensitivity of Distributionally Robust BSDEs}\label{subsec:sensiBSDE}

Let $r \geq 0$. Given a generator $f:\Omega\times[0,T]\times\R\times\R^d\to\R$ and a process $\beta \in \L^\infty$, define
\vspace{-1.5mm}
\begin{equation}\label{eqdef:generatorLinftydeviation}
f^\beta_t (y, z) := f_t (y, z) - \beta_t \cdot z,
\,\, \text{and} \,\ 
F^{r}_t(y, z) := f_t (y, z) + r \vert z \vert.
\vspace{-1.5mm}
\end{equation}
Provided that ${\rm BSDE}(f^\beta,\xi)$ is well posed, with solution $(Y^\beta,Z^\beta)$, we consider the distributionally robust BSDE problem
\vspace{-1.5mm}
\begin{equation}\label{eqdef:sensiLinfty}
  V_{f, \infty}(r):=\sup_{\beta\in\Bc^\infty(r)}Y^\beta_0,
\,\, \text{
where }
\Bc^\infty(r):=\big\{\beta\in\L^\infty(\R^d): \Vert \beta \Vert_\infty\leq r \big\}
.
\vspace{-1.5mm}\end{equation}

\begin{Assumption}\label{ass:controlLinf générales}
\leavevmode
The generator $f:\Omega\times[0,T]\times\R\times\R^d\to\R$ is jointly progressively measurable, satisfies $(f_t(0,0))_{t\in[0,T]}\in\H^{2}$ and fulfills the following conditions:
\begin{romanenumerate}
  \item \label{controlLinfass:f Lip uniformément} $f$ is Lipschitz in $(y, z)$ uniformly in $t$ and $\omega$.
  \item \label{controlLinfass:f C1} For all $t\in[0,T]$, the map $(y, z) \mapsto f_t (y, z)$ is continuously differentiable.
\end{romanenumerate}
\end{Assumption}

\begin{Theorem}\label{thm:V'(0) Linf gen} 
Let $f$ satisfy Assumption \ref{ass:controlLinf générales} and let $\xi \in \L^2$.   
Then,
\begin{romanenumerate}
  \item For all $ r \geq 0$, the ${\rm BSDE}(F^r, \xi)$ admits a unique solution $(Y^r, Z^r)\in\S^2\times\H^2(\R^d)$, and $V_{f, \infty}(r) = Y^r_0$.
  \item The map $ r\mapsto (Y^r, Z^r) \in \S^2 \times \H^2$ is differentiable at $r = 0$, with
$
  (U, V) := \partial_r (Y^r, Z^r) \vert_{r =0},
$
  where $(U,V)$ is the unique solution to the linear ${\rm BSDE}(g^0, 0)$, with
\begin{equation}\label{eqdef:generatorg^0}
g^0_t(u, v) := \vert Z_t \vert + \partial_yf(Y_t,Z_t) u + \partial_zf(Y_t,Z_t)\cdot v, 
\text{where } (Y, Z) := (Y^0, Z^0).
\end{equation}
In particular, letting 
$\Gamma_t^\cdot := \Ec ( \int_t^\cdot \partial_yf(Y_u,Z_u) \mathrm{d}u 
-\int_t^\cdot \partial_z f(Y_u,Z_u) \cdot \mathrm{d}X_u )
,$
\begin{equation}\label{eqdef Gamma t s}
\partial_{r} Y^0_t 
=
\E_t
\Big[ \int_t^T \Gamma_t^s \vert Z_s \vert \mathrm{d}s \Big].
\end{equation}
\end{romanenumerate}
\end{Theorem}

\begin{Remark}\label{rem:conseqassumption}
{\rm
In the problem \eqref{eqdef:sensiLinfty}, one may replace the supremum by an infimum. 
In that case, let $(\hat{Y}, \hat{Z})$ be the unique solution of ${\rm BSDE}(\hat{f},-\xi)$, where
$
\hat{f}_t(y,z):=-f_t(-y,-z),
$
and let $(\hat{Y}^r,\hat{Z}^r)$ be the unique solution of ${\rm BSDE}(\hat{F}^r,-\xi)$, where
$
\hat{F}^r_t(y,z):=\hat{f}_t(y,z)-r|z|.
$
Defining $\hat{\Gamma}$ analogously to \eqref{eqdef Gamma t s}, but with $(\hat{Y},\hat{Z},\hat{f})$ in place of $(Y,Z,f)$, we obtain the analogous representation $
\partial_r \hat{Y}^r_t
=
-\E_t
[
\int_t^T \hat{\Gamma}_t^s |\hat{Z}_s|\,\mathrm{d}s
].
$
}
\end{Remark}
\subsection{Sensitivity of Distributionally Robust Optimal Control}\label{subsect:opticontrol}
We next study the $\L^\infty$ and $\L^2$ sensitivities of optimal control problems. 
The results of this section build on those of the previous one, since it is well known that control problems are intrinsically connected to backward SDEs.
Let $ A \subset \R^d$ and denote by 
\vspace{-2mm}
\begin{equation}\label{eqdef:Avaluedprogmeas}
\text{$\Ac$ the class of progressively measurable processes taking values in $A$.}
\vspace{-2mm}
\end{equation}
Let $ k, l, \lm : \Omega \times [0, T] \times A \rightarrow \R \times \R \times \R^d $.
In the following, for $0 \leq s \leq T $ and $\alpha \in \Ac$, we define 
\vspace{-5.5mm}
\begin{equation}\label{eqdef:constant problems alpha}
\begin{aligned}
l_s^\alpha &:= l(s,\alpha_s),\ 
\lm_s^\alpha := \lm(s,\alpha_s),\ 
k_s^\alpha := k(s,\alpha_s),\ 
\Kc_t^\alpha := e^{-\int_0^t k_s^\alpha \mathrm{d}s},\\
J(\alpha,\beta) &:= \E^{\lm^\alpha+\beta}\Big[\Kc_T^\alpha\xi_T+\int_0^T \Kc_t^\alpha l_t^\alpha \mathrm{d}t\Big],\qquad
J^0(\alpha) := J(\alpha,0).
\end{aligned}
\end{equation}
\vspace{-4.5mm}

\noindent
The definition of $J^0(\alpha)$ is motivated by the fact that many stochastic control problems can be written in the form
$
\inf_{\alpha \in \Ac} J^0(\alpha)
.
$
We define the distributionally robust optimal control (DROC) problem with $\L^\infty$ drift deviation by
\vspace{-2mm}
$$
\overline{V}_{\infty}(r)
:=
\inf_{\alpha\in\Ac}
\sup_{\beta\in\Bc^\infty(r)} J ( \alpha, \beta ) 
,
\vspace{-2mm}
$$ 
where $\Bc^\infty (r)$ is defined in \eqref{eqdef:sensiLinfty}. This is a distributionally robust optimization problem in which the underlying probability measure is regarded as uncertain. At the reference point $\beta =0$, the optimal control problem is formulated under the Wiener measure. 
We consider deviations that are absolutely continuous with respect to the Wiener measure, with density given by the Girsanov transform. 
By reversing the order of the infimum and the supremum, we similarly define
\vspace{-2.5mm}
$$
\underline{V}_{\infty}(r)
:=
\sup_{\beta\in\Bc^\infty(r)}
\inf_{\alpha\in\Ac}
J ( \alpha, \beta ) 
.
\vspace{-4.5mm}
$$ 
In the $\L^2$ case, let
\vspace{-3.5mm}
\begin{equation*}
\Bc^2_{\alpha}(r):=
\Big\{\beta\in\Hb(\R^d):\E^{\lm^\alpha+\beta}
\Big[\int_0^T \Kc_t^\alpha \vert\beta_t\vert^2\mathrm{d}t\Big]\leq r^2\Big\}
\vspace{-3mm}
\end{equation*}
and the corresponding DROC problem is given by
\vspace{-2mm}
\begin{equation}\label{eqdef:L2 sensitivity control}
  \overline{V}_2(r) = \inf_{\alpha\in\Ac}\sup_{\beta\in\Bc^2_{\alpha}(r)}
J ( \alpha, \beta ) .
\vspace{-1mm}
\end{equation}

\begin{Remark}\label{rem:alpha dependance deviation}
{\rm 
 Here, the deviation set $\Bc^2_{\alpha}$ depends on $\alpha$, and therefore the infimum and supremum cannot be interchanged. 
This dependence arises for technical reasons: it allows us to use dynamic programming, although the criterion itself would permit such an inversion if it were independent of $\alpha$. 
We choose this criterion because, in the case $k = 0$, it can be interpreted as an entropy criterion between the reference measure $\P^{\lm^\alpha}$ and the perturbed measure $\P^{\lm^\alpha + \beta}$.
}
\end{Remark}
\vspace{-2mm}
We shall analyze the sensitivity of these distributionally robust optimal control problems by means of backward SDEs. 
We first introduce the generator corresponding to the problem with $r=0$.
\vspace{-1mm}
\begin{equation}\label{eqref:generatorcontrolcase}
f_t(y,z):= \displaystyle\essinf_{a \in A } \{ l^a_t-k^a_ty-\lm^a_t\cdot z \}
\vspace{-1mm}
.\end{equation}

\begin{Assumption}\label{ass:controlLinf cas controle générales}
The maps $l, k, \lm$ are $\F \times \Bc_A$-progressively measurable and satisfy
\leavevmode
\begin{romanenumerate}
  \item \label{controlLinf_control case ass:f=inf falpha}
  For all $t\in[0,T]$,
  $f_t(y,z)$ is continuously differentiable in $(y,z)$.
  \item \label{controlLinf_control case ass:k,lm bornés par L} Both $k$ and $\lm$ are bounded, and $l^\alpha\in\H^{2}$.
  \item \label{existence argmin} There exists a progressively measurable map $\alpha^* : \Omega \times [0, T] \times \R \times \R^d \rightarrow A$ such that 
  \vspace{-2mm}
  $$\arg \displaystyle\essinf_{a \in A}f^a_t(y,z)=\{\alpha^*_t(y,z)\}
  \,\, \text{for all $t \in [0, T]$ and $(y, z) \in \R \times \R^d$}
  .
  \vspace{-3mm}$$
\end{romanenumerate}
\end{Assumption}

If Assumption \ref{ass:controlLinf cas controle générales} is satisfied, then the Hamiltonian $f_t(y,z):=\essinf_{a \in A}f^a_t(y,z)$ satisfies Assumption \ref{ass:controlLinf générales}.
Assumption \ref{ass:controlLinf cas controle générales} is mainly technical and is introduced to ensure the stability estimates required for the $\L^\infty$ sensitivity analysis. 
To satisfy \ref{controlLinf_control case ass:f=inf falpha}--\ref{existence argmin} of Assumption \ref{ass:controlLinf cas controle générales}, one may for instance choose $l$ to be strictly convex in the control variable (for example, quadratic), and $k$ and $\lm$ to be linear, with a bounded convex control set $A$.
We emphasize that no specific structural assumption is imposed on the set $A$. What really matters is the Hamiltonian $f$ associated with the control problem. Sensitivity analysis of $\L^2$ DROC problem requires stronger assumptions.

\begin{Assumption}\label{ass:controlL2 générales}
$l$ is bounded, and the gradient of $f_t(y,z):=\essinf_{a\in A} f^a_t(y,z)$ is $L$-Lipschitz with respect to $(y,z)$, uniformly in $t$ and $\omega$.
\end{Assumption}
\vspace{-1mm}

\noindent Since \citeauthor{el1997backward} \cite{el1997backward} it is known that the optimal control problem 
\vspace{-1mm}
\begin{equation}\label{eqdef:eqdef:opt_cont_pb}
\inf_{\alpha\in\Ac} J^0 (\alpha)
\vspace{-1mm}
\end{equation}
is related to ${\rm BSDE}(f, \xi)$. Namely, ${\rm BSDE}(f, \xi)$ is the non-Markovian counterpart of the Hamilton-Jacobi-Bellman equation. Provided it exists, let $(Y,Z)\in\S^2\times\H^2(\R^d)$ denote the unique solution of ${\rm BSDE}(f , \xi)$. Define 
\vspace{-3mm}
\begin{equation}\label{eqdef:kopt_and_lambdaopt}
\lm^*:=-\partial_z f(Y,Z)
\, \, 
,
\,\, 
k^*:=-\partial_yf(Y,Z)
\, \, 
\text{and} 
\,\, 
\Kc^*_t := e^{- \int_0^t k^*_s \mathrm{d}s}
\vspace{-2mm}
.
\end{equation}

\begin{Theorem}\label{thm:V'(0) Linf et L2} 
\leavevmode
Let Assumption \ref{ass:controlLinf cas controle générales} hold. 
Then, for all $\xi \in \L^2$, ${\rm BSDE}(f,\xi)$ is well defined and admits a unique solution $(Y,Z)\in\S^2\times\H^2(\R^d)$. 
\begin{romanenumerate}
\item\label{thm:V'(0) Linf control} 
Furthermore, for all $r \geq 0$, one has
$
\underline{V}_\infty(r) = \overline{V}_\infty(r) =: V_\infty(r),
$
and the map $V_\infty$ is differentiable at $r=0$, with
\vspace{-4mm}
$$
V_\infty'(0) = \Vert Z \Vert_{\L^1_T(\P^{\lm^*},k^*)}
:= \E^{\lm^*} \Big[ \int_0^T \Kc_s^* \vert Z_s \vert \mathrm{d}s \Big].
\vspace{-3.5mm}
$$
Assume in addition that, for all $t \in [0,T]$, the feedback map $\alpha_t^*(y,z)$ introduced in Assumption 
\ref{ass:controlLinf cas controle générales}-\ref{existence argmin}
is continuously differentiable with respect to $(y,z)$, and that its derivatives $\partial_y \alpha^*$ and $\partial_z \alpha^*$ are uniformly bounded.
Let $r>0$, and let $(Y^r,Z^r)$ denote the solution to ${\rm BSDE}(F^r,\xi)$ (which is well posed) where $F^r_t(y, z) = f_t (y, z) + r \vert z \vert$. Then the associated feedback control
$
\alpha^r := \alpha^*(Y^r,Z^r)
$
is optimal for $\overline{V}_\infty(r)$ and admits the first-order expansion
\vspace{-2mm}
\begin{equation}\label{first order expansion DROcontrol}
\alpha^r
=
\alpha^*(Y,Z)
+
r\,\partial_{y} \alpha^*(Y,Z)\, U
+
r\,\partial_{z} \alpha^*(Y,Z)\cdot V
+
o_{\H^2}(r),
\vspace{-2mm}
\end{equation}
where $(U,V)$ is the solution to the linear ${\rm BSDE}(g^0,0)$ with generator
$
g_t^0(u,v) := k_t^* u + \lm_t^* \cdot v + |Z_t|, 
$
and $o_{\H^2}(r)$ denotes a remainder term such that
$
\|o_{\H^2}(r)\|_{\H^2} = o(r)
\quad \text{as } r \to 0.
$
  \item\label{thm:V'(0) L2} If in addition $\xi \in \L^\infty $ and Assumption \ref{ass:controlL2 générales} holds, then $\overline{V}_2$ is differentiable at $0$ and,
  \vspace{-2.5mm}
  \begin{equation*}
    \overline{V}_2'(0)= \Vert Z \Vert_{\L^2_T(\P^{\lm^*},k^*)} 
    := \E^{\lm^*} \Big[ \int_0^T \Kc_s^* \vert Z_s \vert^2 \mathrm{d}s \Big]^{1/2}.
    \vspace{-3mm}
  \end{equation*}
\end{romanenumerate}
\end{Theorem}

\begin{Remark}\label{rem:linkthm}
{\rm 
\noindent $\bullet$ From the BSDE viewpoint, the first point of the previous theorem clearly appears as a special case of Theorem \ref{thm:V'(0) Linf gen}. 
However, in the second case, namely the $\L^2$ sensitivity setting, as noted in Remark \ref{rem:alpha dependance deviation}, we were not able to define a meaningful deviation set for a general BSDE.

\noindent $\bullet$ In the $\L^\infty$ case, we derived a first-order expansion for the DROC. 
However, for technical reasons, a similar result appears difficult to obtain in the $\L^2$ deviation case; see Remark \ref{rem:foealphaL2}.
}
\end{Remark}
\vspace{-2mm}

Following a referee's suggestion, we also establish an envelope theorem, in the spirit of Theorem $2.14$ of \citeauthor{bartl2023sensitivity} \cite{bartl2023sensitivity}, 
where the authors obtain sharper bounds (of order $O(r^2)$) for the difference between the DROC problem and the standard robust problem evaluated at the optimal control.
To this end, we introduce the two functions
\vspace{-2mm}
$$
V^{\star}_{\infty}(r) := 
\sup_{\beta\in\Bc^\infty(r)}
J(\bar{\alpha}^*, \beta),
\,\, \text{and} \,\,
  V^{\star}_{2}(r) := 
\sup_{\beta\in\Bc^2_{\bar{\alpha}^*}(r)}
J(\bar{\alpha}^*, \beta), \text{where }
\bar{\alpha}^* := \alpha^*(Y,Z),
\vspace{-4mm}
$$
$\alpha^*$ denoting the feedback map introduced in Assumption \ref{ass:controlLinf cas controle générales}, \ref{existence argmin}, and $(Y,Z)$ is the solution of ${\rm BSDE}(f,\xi)$.
In our framework, this envelope theorem is a direct consequence of Theorem \ref{thm:V'(0) Linf et L2}, applied to the uncontrolled case.
\vspace{-2mm}

\begin{Corollary}\label{Envelop Theorem}
\leavevmode
Let Assumption \ref{ass:controlLinf cas controle générales} hold.
\begin{romanenumerate}
\item
For all $\xi \in \L^2$ and all $r \geq 0$, one has
$
V_\infty(r)
=
V^{\star}_{\infty}(r)
+ o(r).
$
  \item
  Under the additional assumptions that $\xi \in \L^\infty$ and Assumption \ref{ass:controlL2 générales} holds, one has
  $
  \overline{V}_2(r)
  =
  V^{\star}_{2}(r)
  + o(r).$
  \vspace{-2mm}
\end{romanenumerate}
\vspace{-3mm}
\end{Corollary}

\subsection{Sensitivity of Distributionally Robust Mixed Problems}
We finally study the $\L^\infty$ and $\L^2$ sensitivities of mixed optimal control/stopping problems:
\vspace{-3mm}
\begin{equation}\label{eqdef:mixedproblem}
\inf_{\tau\in\Tc}
\inf_{\alpha\in\Ac}
J_{\tau}^0 ( \alpha)
,
\vspace{-4mm}
\end{equation}
where $J_\tau$ is defined by
\vspace{-2.5mm}
\begin{equation}\label{eqdef:constant problems alpha DROS}
J_\tau(\alpha,\beta) := \E^{\lm^\alpha+\beta}\Big[\Kc_\tau^\alpha\xi_\tau+\int_0^\tau \Kc_t^\alpha l_t^\alpha \mathrm{d}t\Big],\qquad
J_\tau^0(\alpha) := J_\tau(\alpha,0).
\vspace{-2mm}
\end{equation}
Here $\lm^\alpha,\Kc^\alpha, l^\alpha$ are defined in \eqref{eqdef:constant problems alpha}, $\Ac$ is defined in \eqref{eqdef:Avaluedprogmeas}, $(\xi_t)_{t\in[0,T]}\in\H^0$ satisfies suitable integrability conditions, and $\Tc$ is defined in \eqref{eqdef:Stopping time}. 
One may also consider analogous problems with maximization over the control $\alpha$ instead of the minimization in \eqref{eqdef:mixedproblem}. 
We do not pursue this extension, since it can be treated by the same line of arguments.

We first define the distributionally robust mixed optimal control/stopping problems (DROCS) with $\L^\infty$ deviations:
\vspace{-2mm}
\begin{equation*}
    \underline{\V}_\infty(r) 
    :=\sup_{\beta\in\Bc^\infty(r)}\inf_{\tau\in\Tc}\inf_{\alpha\in\Ac} J_{\tau} (\alpha, \beta) 
\,\, \text{and} \,\,
    \overline{\V}_\infty(r) 
    :=\inf_{\tau\in\Tc}\inf_{\alpha\in\Ac}\sup_{\beta\in\Bc^\infty(r)}
J_{\tau} (\alpha, \beta) 
,
\vspace{-2mm}
\end{equation*}
and the DROCS with $\L^2$ deviations:
\vspace{-2mm}
\begin{equation*}
     \overline{\V}_2(r) :=\inf_{\tau\in\Tc}\inf_{\alpha\in\Ac}\sup_{\beta\in\Bc^2_{\tau, \alpha}(r)}
J_{\tau} (\alpha, \beta)
.
\vspace{-2mm}
\end{equation*}
Here
$
\Bc^2_{\tau, \alpha}(r) :=  \{\beta\in\Hb(\R^d):\E^{\lm^{\alpha} + \beta} [\int_0^\tau \Kc^{\alpha}_t \vert \beta_t \vert ^2\mathrm{d}t ]\leq r^2 \}
$
for all $\tau \in \Tc$ and $\alpha \in \Ac$.
Similarly to the optimal control setting of Subsection \ref{subsect:opticontrol}, we do not consider the corresponding $\sup\inf$ problem, since the set of admissible deviations for $\beta$ depends on the pair $(\tau, \alpha)$. 

We recall the following regularity assumption on the obstacle process $(\xi_t)_{t\in[0,T]}$, which is standard in optimal stopping problems (see, \textit{e.g.}, Assumptions 6 and 7 of \citeauthor{el1997reflected} \cite{el1997reflected}). 
This assumption ensures the existence of an optimal stopping time and will be useful to establish regularity properties of the perturbed stopping time in the DROCS problem.

\begin{Assumption}\label{ass:stop}
The process $(\xi_t)_{t\in[0,T]}\in\H^0$ is continuous a.s., on $[0,T)$ and $\xi_{T^-}:=\lim_{t\to T^-}\xi_t\geq\xi_T$.
\end{Assumption}
\vspace{-1mm}

Consider the generator $f$ defined in \eqref{eqref:generatorcontrolcase}. It is well known, by Theorem $5.10$ of \citeauthor{el2006non} \cite{el2006non}, that the mixed optimal control/stopping problem \eqref{eqdef:mixedproblem} is related to ${\rm RBSDE}(f , \xi)$, which is the non-Markovian counterpart of the free-boundary PDE associated with the mixed optimal control/stopping problem. 
Under the conditions of Theorem \ref{thm:sensimixedproblem} below, ${\rm RBSDE}(f , \xi)$ admits a unique solution $( Y, Z, K) \in\S^2\times\H^2(\R^d) \times \mathbb{A}^2$ and induces the optimal stopping rule and the optimal controlled coefficients
\vspace{-2mm}
\begin{equation}\label{eqdef:tau_opt_mixed}
\tilde{\tau}:=\inf\{t\geq0 \!:\! Y_t \!= \! \xi_t\},
\,\,
\lm^* :=-\partial_z f(Y,Z)
\,
,
\, 
k^*:=-\partial_yf(Y,Z)
\, 
\text{and} 
\,
\Kc^*_t := e^{- \int_0^t k^*_s \mathrm{d}s}
.
\vspace{-1mm}\end{equation}

\begin{Theorem}\label{thm:sensimixedproblem}
Let Assumption \ref{ass:controlLinf cas controle générales} hold, and let $\xi$ satisfy Assumption \ref{ass:stop}. Assume moreover that $\E [\sup_{0\leq t\leq T} \vert \xi_t \vert ^2 ]<+\infty$. Then ${\rm RBSDE}(f, \xi)$ is well posed and
\begin{romanenumerate}
  \item\label{thm:Vup=Vdown mixed inf} for all $r\geq0$, $\overline{\V}_\infty(r)=\underline{\V}_\infty(r) =: \V_\infty(r)$, and $\V_\infty$ is differentiable at $r=0$, with
  \vspace{-2mm}
  \begin{equation*}
     \V_\infty'(0) = 
     \Vert Z \Vert_{ 
     \L^1_{\tilde{\tau}} ( \P^{\lm^*},k^*)
     } 
     := \E^{\lm^*} \Big[ \int_0^{\tilde{\tau }} \Kc_s^* \vert Z_s \vert \mathrm{d}s \Big].
     \vspace{-4mm}
  \end{equation*}
  \item\label{thm:V'(0) mixed L2} If in addition $(\xi_t)_{t\in[0,T]}\in\L^\infty$ and Assumption \ref{ass:controlL2 générales} holds, then $\overline{\V}_2$ is differentiable at $r=0$, and
  \vspace{-3mm}
  \begin{equation*}
    \overline{\V}_2'(0)=
    \Vert Z \Vert_{
    \L^2_{\tilde{\tau}} ( \P^{\lm^* },k^* ) }
    := 
    \E^{\lm^*} 
    \Big[ \int_0^{\tilde{\tau}} \Kc_s^* \vert Z_s \vert^2 \mathrm{d}s \Big]^{1/2}.
    \vspace{-3mm}
  \end{equation*}
\end{romanenumerate}
\end{Theorem}

\begin{Remark}
{\rm 
\noindent $\bullet$ We thank the referee for pointing out that the sensitivity formula in point \ref{thm:V'(0) Linf control} of Theorem \ref{thm:V'(0) Linf et L2} can be viewed as a special case of Theorem \ref{thm:sensimixedproblem}. 
Indeed, by setting $\xi_t = +\infty$ for $t<T$ and $\xi_T = \xi$, one formally recovers Theorem \ref{thm:V'(0) Linf control} from Theorem \ref{thm:sensimixedproblem}. 
However, due to the technical nature of the proofs, we prefer to keep the current presentation, with a first case involving only control, and a second part where optimal stopping is added.

\noindent $\bullet$ A convenient special case of Theorem \ref{thm:sensimixedproblem} is obtained by removing the control, \textit{i.e.}, by taking $f\equiv 0$ (equivalently $A=\{0\}$, so that $\alpha\equiv 0$, $\Kc^\alpha\equiv 1$, and $\lm^\alpha\equiv 0$). 
This yields the distributionally robust optimal stopping problem
$
\inf_{\tau\in\Tc}\sup_{\beta\in\Bc^\infty(r)}\E^{\beta}[\xi_\tau]
$
(and its $\L^2$ analogue with $\beta\in\Bc^2_{\tau, 0}(r)$). 
Under Assumption \ref{ass:stop} and $\E[\sup_{0\le t\le T}|\xi_t|^2]<\infty$, one has $\underline{\V}_\infty(r)=\overline{\V}_\infty(r)=:\V_\infty(r)$ and $\V_\infty$ is differentiable at $r=0$, with
$
\V_\infty'(0)=\Vert \hat{Z}\Vert_{\L^1_{\hat{\tau}}}
$,
where $(\hat{Y},\hat{Z},\hat{K})$ solves ${\rm RBSDE}(0,\xi)$ and $\hat{\tau}:=\inf\{t\ge0:\hat{Y}_t=\xi_t\}$.
If in addition $(\xi_t)_{t\in[0,T]}\in\L^\infty$, then $\overline{\V}_2$ is differentiable at $r=0$ and
$
\overline{\V}_2'(0)=\Vert \hat{Z}\Vert_{\L^2_{\hat{\tau}}}
$.
}
\vspace{-1mm}
\end{Remark}

\begin{Corollary}\label{Envelop Theorem Mixed Pb}
\leavevmode
Let Assumption \ref{ass:controlLinf cas controle générales} hold.
\vspace{-2mm}
\begin{romanenumerate}
  \item\label{thm:Env thm Linf control} 
  Then, for $\xi$ satisfying Assumption \ref{ass:stop}, one has
$
  \V_\infty(r)
  =
  \V^{\star}_{\infty}(r)
  + o(r)$ where
  $
  \V^{\star}_{\infty}(r)
:=
\sup_{\beta\in\Bc^\infty(r)}
J_{\tl{\tau}}(\alpha^*, \beta)$
\vspace{-1mm}

  \item\label{thm:Env thm L2 control}
  Under the additional assumptions that $\xi \in \L^\infty$ and Assumption \ref{ass:controlL2 générales} holds, one has $
  \overline{\V}_2(r)
  =
  \V^{\star}_{2}(r)
  + o(r)
  $
where
  $\V^{\star}_{2}(r)
:= 
\sup_{\beta\in\Bc^2_{\tl{\tau},\alpha^*} (r)}
J_{\tl{\tau}}(\alpha^*, \beta)
  .$
\end{romanenumerate}
\end{Corollary}
\vspace{-5mm}

\section{Literature Comparison and Numerical Illustrations}
\vspace{-1mm}

\subsection{Comparison with Existing Sensitivity Results}

We compare our results to those of \citeauthor{bartl2024numerical} \cite{bartl2024numerical}, \cite{bartl2023sensitivity}, and \citeauthor{jiang2024sensitivity} \cite{jiang2024sensitivity}.
\vspace{-4mm}

\paragraph{Which sensitivity are we considering?}
In \citeauthor{bartl2024numerical} \cite{bartl2024numerical}, \cite{bartl2023sensitivity}, the authors allow for deviations in both drift and volatility. 
By contrast, \citeauthor{jiang2024sensitivity} \cite{jiang2024sensitivity} focuses on martingale deviations, which can be viewed as volatility uncertainty (see also their discussion in Section~4.3 for extensions to the drift case). 
In this paper, we restrict attention to drift deviations only. Our approach is based on a BSDE representation of the DRO problem; incorporating volatility uncertainty would naturally lead to a $2$BSDE framework, which we leave for future work.
As a consequence, the sensitivities in Theorems \ref{thm:V'(0) Linf et L2} and \ref{thm:sensimixedproblem} are not directly comparable to those of \citeauthor{jiang2024sensitivity} \cite{jiang2024sensitivity}. 
We therefore compare our results with \citeauthor{bartl2024numerical} \cite{bartl2024numerical} and \citeauthor{bartl2023sensitivity} \cite{bartl2023sensitivity} in the absence of volatility uncertainty.
\vspace{-4mm}

\paragraph{Comparison of the $\L^2$ sensitivity in the uncontrolled setting, with the Wiener measure as reference model (i.e., at the reference model with zero drift).}

We recall that, given a reference pair $(b^0,\sigma^0)$ satisfying the usual measurability and integrability assumptions,
\citeauthor{bartl2023sensitivity} \cite{bartl2023sensitivity} considered the following value function
\vspace{-3mm}
\begin{equation}\label{eqdef:aux_pbm}
V(r)
:=
\inf_{H\in \Hc}
\sup_{(b,\sigma)\in \Bc_r}
\E\Big[
\int_0^T \!\!\! 
g
\big(t, x_0+ (H \cdot S^{b,\sigma} )_t, H_t\big)\,\mathrm{d}t
+
c
\big(x_0+(H\cdot S^{b,\sigma})_T, \ell(S^{b,\sigma}_T)\big)
\Big].
\vspace{-2mm}
\end{equation}
Here $T>0$ denotes the time horizon and $x_0\in\R$ the initial wealth.
The control $H=(H_t)_{t\in[0,T]}$ is a predictable trading strategy in the admissible set $\Hc$
(see Assumption~2.5 in \citeauthor{bartl2023sensitivity} \cite{bartl2023sensitivity}).
For each pair $(b,\sigma)$, the price process $S^{b,\sigma}$ is defined by
\vspace{-3mm}
$$
\mathrm{d}S_t^{b,\sigma} = b_t\,\mathrm{d}t + \sigma_t\,\mathrm{d}W_t,
\,\, S_0^{b,\sigma}=0,
\text{ and $(H\cdot S^{b,\sigma})_t := \int_0^t H_s\,\mathrm{d}S_s^{b,\sigma}$.}
\vspace{-3mm}
$$
The ambiguity set $\Bc_r$ is centered at $(b^0,\sigma^0)$, for instance
\vspace{-3mm}
$$
\Bc_r
:=
\Big\{
(b,\sigma):\
\|b-b^0\|_{\mathbb L^2([0,T])}
+
\|\sigma-\sigma^0\|_{\mathbb L^2([0,T])}
\le r
\Big\}.
\vspace{-2mm}
$$
The only case in which \eqref{eqdef:aux_pbm} can be directly compared with our control problem \eqref{eqdef:eqdef:opt_cont_pb} is the fixed-control setting.
Concretely, we take $H^\ast\equiv 0$ and $\sigma^0\equiv 1$, which corresponds in our framework to $k^\alpha\equiv 0$ and to a constant drift $\lm^\alpha$.
For simplicity, we further assume $b^0=\lm^\alpha\equiv 0$, so that the reference model is the Wiener measure, where the process is started at $0$, hence $x_0 =0$.

In the $\L^2$ case, our deviation criterion is not norm-based, unlike that of \citeauthor{bartl2023sensitivity} \cite{bartl2023sensitivity};
as discussed in Remark \ref{rem:alpha dependance deviation}, it is closer in spirit to an entropic criterion.
In the present fixed-control setting, the corresponding value functions can be written as
\vspace{-1mm}
$$
V_{\rm bnp}(r)
\!=\!\!\!\!
\sup_{\|b\|_{\H^2}\le r}
\!\!
\E\Big[
\int_0^T 
\!\!\!\!
g(t,S^{b,1}_t)\,\mathrm{d}t
+
u\big(S^{b,1}_T\big)
\Big],
\, \text{and} \,
V_2(r)
=
\!\!\!\!\!\!
\sup_{\Vert \beta \Vert_{\H^2 ( \P^{\beta} )}\le r^2}
\!\!\!
\E^\beta\Big[\!
\int_0^T \!\!\!\!g(t,S^{\beta,1}_t)\,\mathrm{d}t
+
u\big(S^{\beta,1}_T\big)
\Big].
$$
Applying Theorem~2.13 of \citeauthor{bartl2023sensitivity} \cite{bartl2023sensitivity}, Theorem~\ref{thm:V'(0) Linf et L2}, and letting
 $(Y,Z)$ be the unique solution of
$
\mathrm{d}Y_t=-g(t,W_t)\,\mathrm{d}t+Z_t\,\mathrm{d}W_t
$, $Y_T=u(W_T)
$, we obtain
\vspace{-1mm}
\begin{equation*}
V_{\rm bnp}'(0)
=
\Big(
\E 
\Big[
\int_0^T
\Big|
\E_t \Big[
u'(W_T)
\!+\!\!
\int_t^T \!\! \partial_x g(s,W_s)\,\mathrm{d}s
\Big]
\Big|^2
\mathrm{d}t
\Big]
\Big)^{1/2} \, ,\,
V_2'(0)
=
\Big(
\E \Big[
\int_0^T |Z_t|^2\,\mathrm{d}t
\Big]
\Big)^{1/2},
\vspace{-2mm}
\end{equation*}
Assuming that $g$ and $u$ are deterministic and sufficiently regular, one has $Y_t=\varphi(t,W_t)$ and $Z_t=\partial_x\varphi(t,W_t)$,
where $\varphi$ satisfies 
the terminal condition $\varphi(T,\cdot)=u$ and
\vspace{-1mm}
$$
\partial_t\varphi+\tfrac12\partial_{xx}\varphi+g=0
\vspace{-1mm}
.$$
By the Feynman--Kac formula,
$
Z_t
=
\E_t [
u'(W_T)
+
\int_t^T \partial_x g(s,W_s)\,\mathrm{d}s
]
$, which proves that $V_2'(0) = V_{\rm bnp}'(0)$.
Hence, despite the different deviation criteria, the resulting $\L^2$ sensitivities coincide in this setting.
The same conclusion holds when $b^0$ and $\lm^\alpha$ are constant, and similarly in higher dimension.
However, this argument breaks down when $b^0$ and $\lm^\alpha$ are stochastic processes.
\vspace{-3mm}

\paragraph{Comparison of the $\L^\infty$ sensitivity in the uncontrolled setting, with the Wiener measure as reference model.}
A similar comparison can be performed in the $\L^\infty$ setting considered in 
\citeauthor{bartl2024numerical} \cite{bartl2024numerical}. 
In the reference case $b^0=0$ and $\sigma^0=1$, they study
\vspace{-4mm}
$$
\tilde V_{\rm bnp}(r)
:=
\sup_{\|\beta\|_\infty\le r}
\E \Big[
c
\Big(\int_0^T \beta_s\,\mathrm{d}s + W_T, 
\ell \Big(
\int_0^T \beta_s\,\mathrm{d}s + W_T
\Big)\Big)
\Big],
\vspace{-4mm}
$$
and show that
\vspace{-2mm}
$$
\tilde V_{\rm bnp}'(0)
=
\int_0^T
\E \Big[
\big|\E_t[\partial_x c(W_T,\ell(W_T))]\big|
\Big]\mathrm{d}t.
$$
On the other hand, applying Theorem \ref{thm:V'(0) Linf et L2} yields
\vspace{-1mm}
$$
V_\infty'(0)
=
\tilde V_{\rm bnp}'(0)
\vspace{-2mm}
$$
Hence the two sensitivities coincide. This is consistent with the fact that, in this setting,
$V_\infty(r)=\tilde V(r)$.

\subsection{Numerical Illustrations}

\subsubsection{Hedging Under Different Borrowing-Lending Rates}\label{subsec:borrow_lend}

We illustrate our notion of sensitivity for distributionally robust optimal
control problems through an optimal hedging problem with different borrowing and lending rates.

Let $\xi$ be the payoff at maturity $T$ of a derivative security in a financial market
with borrowing and lending rates satisfying $\gm_{\min}<\gm_{\max}$.
For a given short-rate process $\alpha$, the corresponding pricing functional is given by
\vspace{-4mm}
$$
\mathrm{Price}^{(\alpha)}(\xi)
:=
\E^{\alpha}\!\big[
\Kc_T^\alpha\,\xi
\big],
\qquad
\Kc_t^\alpha
:=
\exp\!\Big(-\int_0^t \alpha_s\,\mathrm{d}s\Big).
\vspace{-3.5mm}
$$
In order to ensure sufficient regularity of the Hamiltonian, we consider a penalized
version of the pricing problem.
Fix $a>0$, let $
\bar \gm := \frac{\gm_{\min}+\gm_{\max}}{2}
$ and define
\vspace{-3mm}
$$
J^a (\alpha, \beta)
:=
\E^{\alpha + \beta}
\!\Big[
\Kc_T^\alpha\,\xi
+
a\,
\int_0^T
\Kc_s^\alpha\,
\big(\alpha_s-\bar \gm\big)^2
\,\mathrm{d}s
\Big].
\vspace{-3mm}
$$
The additional quadratic running cost penalizes deviations of the short rate $\alpha$
from the reference level $\bar\gm$ and ensures the smoothness and convexity properties
required for our analysis.
When the penalization parameter $a$ is small, $\inf_{\alpha \in \Ac} J^a (\alpha, 0)$ can be interpreted as a
regularized approximation of the original pricing problem with free choice of the
short rate.
By contrast, when $a$ is sufficiently large, the optimal control $\alpha$ is forced to
remain close to $\bar\gm$, and the problem asymptotically reduces to the standard pricing
formula associated with the constant short rate $\bar\gm$.

We now introduce the distributionally robust version of the penalized pricing problem.
For $r\ge 0$, define
\vspace{-3mm}
\begin{equation*}
V^a_{\infty}(r)
:=
\inf_{\alpha\in\Ac}
\sup_{\beta \in \Bc^{\infty}(r)}
J^a (\alpha, \beta),
\qquad
V^a_{2}(r)
:=
\inf_{\alpha\in\Ac}
\sup_{\beta \in \Bc^{2}_{\alpha}(r)}
J^a (\alpha, \beta).
\vspace{-3mm}
\end{equation*}
We illustrate these problems with a European lookback option in the Black--Scholes model, whose terminal payoff is
\vspace{-3mm}
$$
\xi = \varphi(M_T-K),
\qquad
M_T := \sup_{0\le t\le T} S_t,
\qquad
S_t = S_0 \exp\Big(\big(\bar\gm-\frac{\sigma^2}{2}\big)t + \sigma W_t\Big),
\vspace{-3mm}
$$
where $(S_t)_{t\in[0,T]}$ is a geometric Brownian motion under the reference measure $\P^0$.
In this case,
$
f(t,y,z)
=
\inf_{\rho\in[\gm_{\min},\gm_{\max}]}
 \{
a\,(\rho-\bar\gm)^2
-
y\,\rho
-
z\,\rho
 \},
$
and $\alpha_t^*(y, z)$ is the corresponding minimizer.
We then define the process $\alpha^* := \alpha_t^*(Y_t, Z_t)$,
where the pair $(Y,Z)$ is the unique solution of the backward stochastic differential
equation ${\rm BSDE}(f,\xi)$.
In this case, a straightforward computation shows that the Hamiltonian admits a Lipschitz continuous gradient.
\vspace{-2mm}

\begin{Proposition}
Assume that $a>0$, $\gm_{\min}<\gm_{\max}$, and that $\varphi$ is bounded and measurable. Then Assumptions~\ref{ass:controlLinf cas controle générales} and~\ref{ass:controlL2 générales} are satisfied. Consequently, Theorem~\ref{thm:V'(0) Linf et L2} applies to both the $\L^\infty$ and $\L^2$ sensitivities.
\end{Proposition}
\noindent\textit{Proof.}
The control set is compact, while $k^\rho=\lm^\rho=\rho$ and $l^\rho=a(\rho-\bar\gm)^2$ are bounded. Strict convexity in $\rho$ yields the unique measurable minimizer
$
\alpha^*(y,z)=\Pi_{[\gm_{\min},\gm_{\max}]}
(\bar\gm+(y+z)/(2a)).
$
Thus, $f$ is continuously differentiable and, since the projection is Lipschitz, so is $\nabla f$. Finally, the boundedness of $\varphi$ implies that $\xi\in\L^\infty$. \hfill$\square$

For the numerical experiment, $\varphi$ is chosen as a smooth and bounded approximation of the map $x\mapsto x^+$.
The associated backward stochastic differential equation ${\rm BSDE}(f,\xi)$ is solved
numerically using a deep BSDE scheme, in which the process $Z$ is parameterized by a
neural network.
We refer to \citeauthor{han2018solving} \cite{han2018solving},
\citeauthor{hure2020deep} \cite{hure2020deep}, and
\citeauthor{han2017deep} \cite{han2017deep} for details on this methodology.
Since the payoff depends on the running maximum of the underlying asset, the process
$(M_t)_{t\in[0,T]}$ is included as an additional state variable in the numerical
implementation.
For the numerical experiments, we compute the sensitivity for several values of the
strike $K$. Parameters are fixed as follows:
\vspace{-2mm}
$$
S_0 = 10,
\quad
T = 0.5,
\quad
\sigma = 0.2,
\quad
\gm_{\min} = 0.02,
\quad
\gm_{\max} = 0.08,
\quad
\bar\gm = 0.04,
\quad
a = 0.1.
$$
\vspace{-10mm}
\begin{figure}[H]
\begin{centering}
\includegraphics[scale = 0.14]{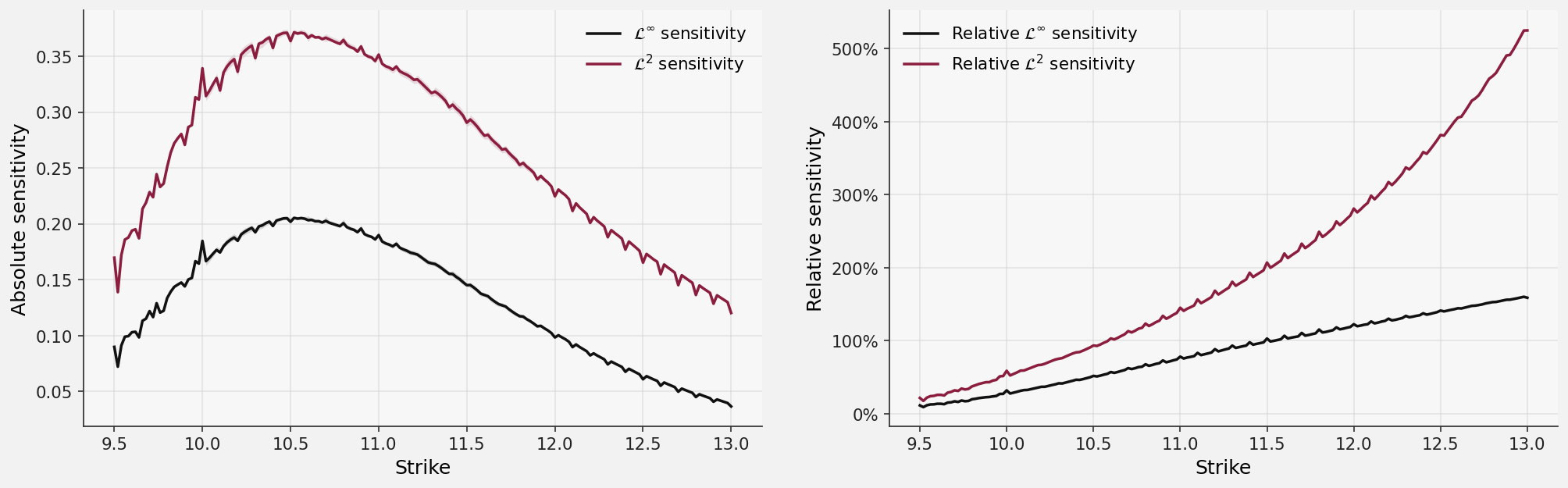}
\caption{\it \footnotesize Sensitivities.}
\end{centering}
\vspace{-3mm}
\end{figure}

\noindent We observe that the two sensitivities display very similar patterns. However, once normalized by the option price, the relative sensitivities reveal that model risk remains substantial under both the $\L^2$ and $\L^\infty$ deviation criteria. In particular, for strikes that are deep out of the money, the relative sensitivities exceed $100\%$. This is mainly due to the fact that the option price tends to zero in this regime.
\vspace{-3mm}

\subsubsection{American lookback option problem}

In the context of the same setting as the previous subsection \ref{subsec:borrow_lend}, with the constant short rate $\gamma>0$,
the price of an American lookback option with strike $K$ is given by the optimal stopping problem
\vspace{-3mm}
$$
C
=
\sup_{\tau\in\Tc}
\E
\big[
e^{-\gm\tau}\,
\varphi(M_\tau-K)
\big].
\vspace{-1mm}
$$
We consider the corresponding distributionally robust pricing problems
\vspace{-2mm}
$$
V^{\mathrm{Lb}}_{\infty}(r)
:=
\sup_{\tau\in\Tc}
\inf_{\beta\in\Bc^{\infty}(r)}
\E^{\beta}
\big[
e^{-\gm\tau}
\varphi(M_\tau-K)
\big]
\,\, \text{and} \,\, 
V^{\mathrm{Lb}}_{2}(r)
:=
\sup_{\tau\in\Tc}
\inf_{\beta\in\Bc^{2}_{\tau}(r)}
\E^{\beta}
\big[ 
e^{-\gm\tau}
\varphi(M_\tau-K)
\big].
$$
\vspace{-4mm}

\noindent
In both cases, the stopping time $\tau$ represents the exercise strategy of the option holder, while the process $\beta$ models an adversarial perturbation of the reference pricing measure.
The quantities $V^{\mathrm{Lb}}_{\infty}(r)$ and $V^{\mathrm{Lb}}_{2}(r)$ can therefore be interpreted as robust prices of the American lookback option under drift uncertainty.
We compute the corresponding sensitivities for several values of the strike $K$.
The payoff function $\varphi$ is chosen as a smooth and bounded approximation of the positive part.

\begin{Proposition}
Assume that $\varphi$ is bounded and continuous. Then Assumptions~\ref{ass:controlLinf cas controle générales}, \ref{ass:stop}, and~\ref{ass:controlL2 générales} are satisfied in the uncontrolled case $f\equiv0$. Consequently, Theorem~\ref{thm:sensimixedproblem} applies to both the $\L^\infty$ and $\L^2$ sensitivities.
\end{Proposition}
\noindent\textit{Proof.}
Taking $A=\{0\}$ makes the control assumptions immediate. In the Black--Scholes model, $S$ is an exponential of Brownian motion and therefore has continuous paths; its running maximum $M$ is continuous as well. Hence the obstacle $\xi_t=e^{-\gm t}\varphi(M_t-K)$ is continuous and bounded. In particular, $\xi_{T^-}=\xi_T$, and all the required integrability conditions hold. \hfill$\square$

The associated reflected BSDE is solved numerically using a deep BSDE scheme.
More precisely, the Skorokhod reflection condition is enforced through a penalization approach, by adding a penalty term to the loss function that discourages violations of the constraint $Y_t\le\xi_t$ along the simulated paths. 
\vspace{-2mm}
$$
\text{The parameters are fixed as follows: }
S_0 = 10,
\qquad
T = 1,
\qquad
\sigma = 1.0,
\qquad
\gm = 1.0.
\vspace{-4mm}
$$

\begin{figure}[H]
\begin{centering}
\includegraphics[scale = 0.35]{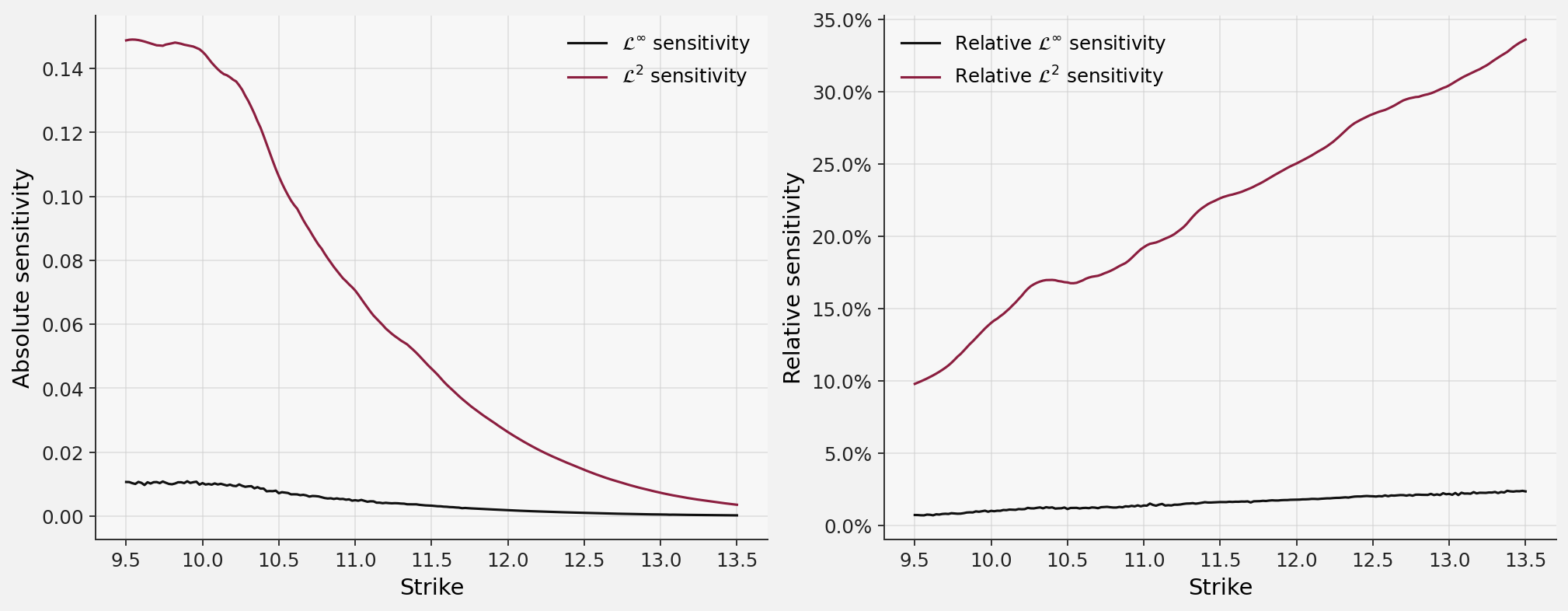}
\caption{\it \footnotesize Sensitivities.}
\end{centering}
\end{figure}
Interestingly, the absolute sensitivity decreases as the strike increases. 
However, when considering relative sensitivities, the decrease in the option price appears to be stronger than the decrease in the sensitivity itself, so that the relative sensitivity increases with the strike. 
We also observe that the American counterpart of the lookback option seems to decay more slowly for deep out-of-the-money strikes than the European counterpart studied previously, as reflected by the slower increase of the relative sensitivity.
\vspace{-3mm}

\section{Proof of the Sensitivities for Control Problems}

\subsection{$\L^\infty$ Sensitivity of Distributionally Robust BSDE}

Let $f:\Omega\times[0,T]\times\R\times\R^d\to\R$. 
We study the differentiability at the origin of the map $V_{f,\infty}$ defined in \eqref{eqdef:sensiLinfty}. 
For $r\geq0$ and $\beta\in\Bc^\infty(r)$, let $F^r$ and $f^\beta$ be defined by \eqref{eqdef:generatorLinftydeviation}.

\smallskip
\noindent \textit{Proof of Theorem \ref{thm:V'(0) Linf gen}.}
\noindent {\bf Step $1$.}
We first prove that
$
V_{f,\infty}(r)=Y_0^r,
$
where $(Y^r,Z^r)$ denotes the unique solution of ${\rm BSDE}(F^r,\xi)$. 
By the standard well-posedness result for BSDEs with Lipschitz generator (Theorem 4.3.1 in \citeauthor{Jianfen_BSDE} \cite{Jianfen_BSDE}), both ${\rm BSDE}(f^\beta,\xi)$ and ${\rm BSDE}(F^r,\xi)$ admit a unique solution in $\S^2\times\H^2(\R^d)$, denoted respectively by $(Y^\beta,Z^\beta)$ and $(Y^r,Z^r)$. 
Observe that
$
F_t^r(Y_t^r,Z_t^r)=f_t(Y_t^r,Z_t^r)+r|Z_t^r|
\geq f_t(Y_t^r,Z_t^r)-\beta_t\cdot Z_t^r
$
for all $\beta\in\Bc^\infty(r)$. By the comparison theorem (Theorem 4.4.1 in \citeauthor{Jianfen_BSDE} \cite{Jianfen_BSDE}), this yields $Y_t^r\geq Y_t^\beta$ for all $t\in[0,T]$, and therefore $Y_0^r\geq V_{f,\infty}(r)$. Conversely, defining
$
\beta_t^r:=-r\,\frac{Z_t^r}{|Z_t^r|}\,\mathbf{1}_{\{Z_t^r\neq0\}},
$
we have $\beta^r\in\Bc^\infty(r)$ and
$
F_t^r(Y_t^r,Z_t^r)=f_t(Y_t^r,Z_t^r)-\beta_t^r\cdot Z_t^r,
$
so that $(Y^r,Z^r)$ solves ${\rm BSDE}(f^{\beta^r},\xi)$ and thus $Y_0^r=Y_0^{\beta^r}\leq V_{f,\infty}(r)$. Hence $V_{f,\infty}(r)=Y_0^r$.

\noindent {\bf Step $2$.}
We now prove that the map
$
r\longmapsto (Y^r,Z^r)\in\S^2\times\H^2(\R^d)
$
is differentiable at $r=0$. 
Let $(Y,Z):=(Y^0,Z^0)$, define
$
U^r:=\frac{Y^r-Y}{r}
$ and $
V^r:=\frac{Z^r-Z}{r}.
$
Then $(U^r,V^r)$ solves ${\rm BSDE}(g^r,0)$, where
\vspace{-3mm}
$$
g_t^r(u,v)
:=
\frac{f_t(Y_t+ru,Z_t+rv)-f_t(Y_t,Z_t)}{r}
+
|Z_t+rv|
\,\, \text{satisfies }
g_t^r(u,v)\xrightarrow[r\to0]{}g_t^0(u,v)
\textit{a.s.},
\vspace{-2mm}
$$
where $g^0$ is defined in \eqref{eqdef:generatorg^0}, and the convergence follows from Assumption \ref{ass:controlLinf générales}-\ref{controlLinfass:f C1}.
Moreover, by Assumption \ref{ass:controlLinf générales}-\ref{controlLinfass:f Lip uniformément}, both $g^r$ and $g^0$ are Lipschitz in $(u,v)$, uniformly in $(t,\omega)$. 
Finally,
$
g_t^r(0,0)=|Z_t|=g_t^0(0,0),
$
and since $Z\in\H^2(\R^d)$, it follows that
$
(g^r(\cdot,\cdot))_{r\geq0}
$ and $
g^0(\cdot,\cdot)
$
satisfy the assumptions of the stability theorem for Lipschitz BSDEs (Theorem 4.4.3 in \citeauthor{Jianfen_BSDE} \cite{Jianfen_BSDE}). 
Therefore,
$
\Vert U^r-U\Vert_{\S^2}\longrightarrow0$ and 
$
\Vert V^r-V\Vert_{\H^2}\longrightarrow0,
$
where $(U,V)$ is the unique solution of ${\rm BSDE}(g^0,0)$.
It remains to identify $U$. Since ${\rm BSDE}(g^0,0)$ is linear, the classical representation formula yields
$
U_t
=
\E_t [\int_t^T \Gamma_t^s |Z_s|\,\mathrm{d}s].$
The process $V$ is then obtained from the martingale representation property.
\ep 
\vspace{-2mm}

\newpage
\subsection{Deriving the $\L^\infty$ Sensitivity of DROC Problems}\label{subsec:proofthemcontrol}

For $r\geq0$, $a \in A$, and $b\in\R^d$ with $|b|\leq r$, recall the generators
\vspace{-2mm}
\begin{equation}\label{eqdef:genLinftyopticont}
\begin{split}
  &f_t(y,z) := \essinf_{a \in A} f_t^a(y,z),
  \qquad
  f_t^a(y,z) := l_t^a - k_t^a y - \lm_t^a \cdot z,
  \\
  &F_t^r(y,z)
  := \sup_{|b|\leq r}\inf_{a \in A} h_t^{a,b}(y,z)
  = f_t(y,z)+r|z|,
  \qquad
  h_t^{a,b}(y,z) := f_t^a(y,z)- b \cdot z.
\end{split}
\end{equation}
  \vspace{-4mm}

\noindent For $\alpha \in \Ac$ and $\beta \in \Bc^{\infty}(r)$, we write, by abuse of notation,
$
h_t^{\alpha,\beta} := h_t^{\alpha_t,\beta_t}
$ and $
f_t^\alpha := f_t^{\alpha_t}.
$

\noindent \textit{Proof of Theorem \ref{thm:V'(0) Linf et L2}, \ref{thm:V'(0) Linf control}.}
\textbf{Step $1$.}
We first prove that
$
\overline{V}_\infty(r)=\underline{V}_\infty(r)=Y_0^r.
$
By Assumption \ref{ass:controlLinf cas controle générales},
\ref{controlLinf_control case ass:f=inf falpha} and
\ref{controlLinf_control case ass:k,lm bornés par L},
the generator $F^r$ is Lipschitz in $(y,z)$.
Hence ${\rm BSDE}(F^r,\xi)$ admits a unique solution
$(Y^r,Z^r)\in\S^2\times\H^2(\R^d)$.
Moreover, by Assumption \ref{ass:controlLinf cas controle générales},
\ref{controlLinf_control case ass:k,lm bornés par L},
for each $\alpha \in \Ac$ and $\beta \in \Bc^\infty(r)$, the linear
${\rm BSDE}(h^{\alpha,\beta},\xi)$ has bounded coefficients and therefore admits a unique solution
$(Y^{\alpha,\beta},Z^{\alpha,\beta})\in\S^2\times\H^2(\R^d)$, with
\vspace{-3mm}
$$
Y_0^{\alpha,\beta}
=
\E^{\lm^\alpha+\beta}\Big[\Kc_T^\alpha\xi+\int_0^T \Kc_t^\alpha l_t^\alpha\,\mathrm{d}t\Big].
\vspace{-3mm}
$$
By Assumption \ref{ass:controlLinf cas controle générales}, \ref{existence argmin}, there exists a progressively measurable map $\alpha^*$ such that, for all $(t,\omega,y,z)$, $\alpha_t^*(y,z)$ is a minimizer of $a\mapsto f_t^a(y,z)$.
Define
$
\hat{\alpha}:=\alpha^*(Y^r,Z^r)\in\Ac$, $
\hat{\beta}:=-r\frac{Z^r}{|Z^r|}\mathbf{1}_{\{Z^r\neq0\}}\in\Bc^\infty(r).
$
Then
$
h_t^{\hat{\alpha},\hat{\beta}}(Y_t^r,Z_t^r)=F_t^r(Y_t^r,Z_t^r).
$
By uniqueness, it follows that
$
Y_t^r=Y_t^{\hat{\alpha},\hat{\beta}}
$ and $
Z_t^r=Z_t^{\hat{\alpha},\hat{\beta}}
 \mathrm{d}t\times \mathrm{d}\P^0\text{-a.e.}
$.
Moreover, for all $\alpha\in\Ac$ and $\beta\in\Bc^\infty(r)$,
$$
h_t^{\hat{\alpha},\beta}(Y_t^{\hat{\alpha},\hat{\beta}},Z_t^{\hat{\alpha},\hat{\beta}})
\leq
h_t^{\hat{\alpha},\hat{\beta}}(Y_t^{\hat{\alpha},\hat{\beta}},Z_t^{\hat{\alpha},\hat{\beta}})
\leq
h_t^{\alpha,\hat{\beta}}(Y_t^{\hat{\alpha},\hat{\beta}},Z_t^{\hat{\alpha},\hat{\beta}})
\qquad 
\qquad \mathrm{d}t\times\mathrm{d}\P^0\text{-a.e.}
$$
Applying the comparison theorem for Lipschitz BSDEs
(Theorem 4.4.1 in \citeauthor{Jianfen_BSDE} \cite{Jianfen_BSDE}), we obtain
$$
Y_0^{\hat{\alpha},\beta}\leq Y_0^{\hat{\alpha},\hat{\beta}}\leq Y_0^{\alpha,\hat{\beta}}.
$$
Hence $(\hat{\alpha},\hat{\beta})$ is a saddle point for
$(\alpha,\beta)\mapsto Y_0^{\alpha,\beta}$,
and therefore
$
\overline{V}_\infty(r)=\underline{V}_\infty(r)=Y_0^r.
$

\noindent \textbf{Step $2$.}
We now prove differentiability at the origin.
By Step $1$, one has
$
\overline{V}_\infty(r)=\underline{V}_\infty(r)=Y_0^r.
$
Since $\Ac$, $l$, $k$, and $\lm$ satisfy Assumption \ref{ass:controlLinf cas controle générales}, the generator $f_t(y,z)$ satisfies Assumption \ref{ass:controlLinf générales}. 
Hence the conclusion follows directly from Theorem \ref{thm:V'(0) Linf gen}.
\ep

The first-order expansion \eqref{first order expansion DROcontrol} for the optimal feedback control follows directly from the previous argument together with the proof of Theorem \ref{thm:V'(0) Linf gen}. 

\subsection{Deriving the $\L^2$ Sensitivity of DRO Problems}

We recall that $\overline{V}_2$ is defined in \eqref{eqdef:L2 sensitivity control}. 
By definition of $\Bc_\alpha^2(r)$, we may rewrite $\overline{V}_2(r)$ as
\begin{align*}
\overline{V}_2(r)
&=
\inf_{\alpha\in\Ac}
\sup_{\beta\in\Hb(\R^d)}
\inf_{\gm>0}
\Big\{
J(\alpha,\beta)
+
\gm
\Big(
r^2
-
\E^{\lm^\alpha+\beta}
\Big[
\int_0^T \Kc_t^\alpha |\beta_t|^2\,\mathrm{d}t
\Big]
\Big)
\Big\}.
\end{align*}
After the change of variable $\gm \leftrightarrow \frac{1}{4\gm}$, define
\vspace{-2mm}
$$
\hat{J}_\gm(\alpha,\beta)
:=
J(\alpha,\beta)
-
\frac{1}{4\gm}
\E^{\lm^\alpha+\beta}
\Big[
\int_0^T \Kc_t^\alpha |\beta_t|^2\,\mathrm{d}t
\Big].
\vspace{-3mm}
$$
We then introduce
\vspace{-2mm}
\begin{equation}\label{eqdef:hatJ}
\overline{V}_2(r)
=
\inf_{\alpha\in\Ac}
\sup_{\beta\in\Hb(\R^d)}
\inf_{\gm>0}
\Lc_r^\alpha(\gm,\beta),
\qquad
\Lc_r^\alpha(\gm,\beta)
:=
\hat{J}_\gm(\alpha,\beta)
+
\frac{r^2}{4\gm}.
\vspace{-2mm}
\end{equation}

\noindent Interchanging the infimum over $\gm$ with the infimum over $\alpha$ and the supremum over $\beta$ naturally leads to the auxiliary value function
\begin{equation}\label{eqdef:tl V(r)}
\tl V_2(r)
:=
\inf_{\gm>0}
\Big\{
G(\gm)
+
\frac{r^2}{4\gm}
\Big\},
\qquad
G(\gm)
:=
\inf_{\alpha\in\Ac}
\sup_{\beta\in\Hb(\R^d)}
\hat{J}_\gm(\alpha,\beta).
\end{equation}

\noindent
We first analyze the function $G$ through a BSDE representation, extending its definition to $\gm=0$ in order to establish differentiability at the origin. 
We will then use the following strong duality lemma to prove that
$
\tl V_2(r)=\overline{V}_2(r).
$

\begin{Lemma}\label{lem:strong duality}
Let $\Xc$ be a set, and let $\psi:\Xc\to\R$ and $\varphi:\Xc\to\R_+$.
Assume that
\vspace{-2mm}
$$
H(\lm)
:=
\sup_{x \in \Xc}
\big\{
\psi(x)-\lm \varphi(x)
\big\}
\vspace{-2mm}
$$
is bounded by some constant $C$, differentiable on $\R_+^*$, and that there exists a family $(x_\lm)_{\lm>0}$ such that, for all $\lm>0$,
\vspace{-2mm}
\begin{equation}\label{eqde:familyoptim}
x_\lm
\in
\argmax_{x \in \Xc}
\big\{
\psi(x)-\lm \varphi(x)
\big\}
\qquad \text{and} \qquad
H'(\lm)=-\varphi(x_\lm).
\vspace{-2mm}
\end{equation}
Then, for all $r>0$, the following strong duality holds:
\vspace{-3mm}
$$
\sup_{x \in \Xc , \varphi(x)\leq r} \psi(x)
=
\inf_{\lm > 0}
\big\{
H(\lm)+\lm r
\big\}.
$$
\end{Lemma}
\vspace{-3mm}
\noindent The proof of Lemma \ref{lem:strong duality} is deferred to Subsection \ref{subsect:DeterministicLemma}. 
This lemma shows that for a constrained optimization problem of the form
$
\sup_{\varphi(x)\leq r} \psi(x),
$
with $\varphi$ nonnegative and $\psi$ bounded, the validity of an envelope formula for $H$, together with a suitable regularity assumption, is sufficient to ensure strong duality. 
Note that no assumption is imposed on $\Xc$, and no convexity assumption is required on either $\psi$ or $\varphi$. 
We now state the final deterministic lemma, which will yield the sensitivity at the origin of $\overline{V}_2$.

\begin{Lemma}\label{lem:concluding lemma}
Let $g:\R_+\to\R$ be a nondecreasing function differentiable at $\gm=0$.
Define, for $r\geq0$,
\vspace{-3mm}
$$
V(r):=\inf_{\gm>0}\Big\{g(\gm)+\frac{r^2}{\gm}\Big\}.
\vspace{-1mm}
$$
Then $V$ is differentiable at $r=0^+$ and
$
V'(0)=2\sqrt{g'(0)}.
$
\end{Lemma}

\vspace{-1mm}
\noindent The proof of Lemma \ref{lem:concluding lemma} is also deferred to Subsection \ref{subsect:DeterministicLemma}. 

\vspace{-1mm}
\begin{Proposition}\label{prop:diff tl V2}
Let Assumption \ref{ass:controlL2 générales} hold, and let $\xi \in \L^\infty$.
Then the function $G$, defined in \eqref{eqdef:tl V(r)}, admits an extension to $[0,+\infty)$ as a nondecreasing function, differentiable at $0$, and, letting $(Y,Z)$ be the solution of ${\rm BSDE}(f,\xi)$, we have
\vspace{-2mm}
$$
G'(0)=\Vert Z \Vert^2_{\L^2(\P^{\lm^*},k^*)}.
\vspace{-2mm}
$$
\end{Proposition}

\vspace{-1mm}
\begin{Remark}
{\rm
\noindent $\bullet$ As a direct consequence of Proposition \ref{prop:diff tl V2} and Lemma \ref{lem:concluding lemma}, one obtains
\vspace{-2mm}
$$
\tl V_2'(0)=\Vert Z \Vert_{\L^2(\P^{\lm^*},k^*)}.
\vspace{-2mm}
$$

\noindent $\bullet$ Proposition \ref{prop:diff tl V2} only yields differentiability of $G$ at $0$, which is why Lemma \ref{lem:concluding lemma} is needed.
}
\end{Remark}

\vspace{-1mm}
\proof
\textbf{Step 1.}
We show that $G(\gm)=Y^\gm_0$, where $(Y^\gm,Z^\gm)$ is the unique solution of ${\rm BSDE}(\hat{F}^\gm,\xi)$.
We now introduce the generators
\vspace{-2mm}
\begin{align*}
&f_t(y,z)
:=
\essinf_{\alpha \in \Ac}
\big\{
l^\alpha_t-k^\alpha_t y-\lm^\alpha_t\cdot z
\big\}\,\, ,\,\,
\hat{h}^{\gm,\alpha,\beta}_t(y,z)
:=
l^\alpha_t-k^\alpha_t y-(\lm^\alpha_t+\beta_t)\cdot z-\frac{1}{4\gm}|\beta_t|^2,
\\
&\text{and}\,
\hat{F}^\gm_t(y,z)
:=
\sup_{b\in\R^d}
\big\{
f_t(y,z)-b\cdot z-\frac{1}{4\gm}|b|^2
\big\}
=
f_t(y,z)+\gm |z|^2.
\end{align*}
\vspace{-4mm}

\noindent By Assumption \ref{ass:controlL2 générales} and since $\xi\in\L^\infty$, we can apply the well-posedness theorem for quadratic BSDEs (Theorem 7.3.3 in \citeauthor{Jianfen_BSDE} \cite{Jianfen_BSDE}) together with the a priori estimate of Theorem 7.2.1 therein. This yields that ${\rm BSDE}(\hat{F}^\gm,\xi)$ admits a unique solution $(Y^\gm,Z^\gm)\in\S^\infty\times\Hb(\R^d)$.
Moreover, for each $\alpha\in\Ac$ and $\beta\in\Hb(\R^d)$, the generator $\hat{h}^{\gm,\alpha,\beta}$ and the terminal condition $\xi$ satisfy the assumptions of the linear BSDE results in \citeauthor{JacksonLinearBSDE_BMO} \cite{JacksonLinearBSDE_BMO} (combining Proposition 5.8, Proposition 4.7 and Theorem 1.5). Hence, the linear ${\rm BSDE}(\hat{h}^{\gm,\alpha,\beta},\xi)$ has a unique solution $(Y^{\gm,\alpha,\beta},Z^{\gm,\alpha,\beta})\in\S^\infty\times\Hb(\R^d)$.

\noindent Let $(\hat\alpha^\gm,\hat\beta^\gm)$ be the saddle-point controls associated with $\hat{F}^\gm$, so that
\vspace{-2mm}
$$
\hat{h}^{\gm,\hat\alpha^\gm,\hat\beta^\gm}_t(Y^\gm_t,Z^\gm_t)
=
\hat{F}^\gm_t(Y^\gm_t,Z^\gm_t)
\qquad \mathrm{d}t\times\mathrm{d}\P^0\text{-a.e.}
\vspace{-2mm}
$$
By uniqueness of solutions to BSDEs with generator $\hat{h}^{\gm,\hat\alpha^\gm,\hat\beta^\gm}$ and terminal condition $\xi$, we obtain
$
Y^\gm_t = Y^{\gm,\hat\alpha^\gm,\hat\beta^\gm}_t
$
for all $t\in[0,T]$ and
$
Z^\gm_t = Z^{\gm,\hat\alpha^\gm,\hat\beta^\gm}_t
\mathrm{d}t\times\mathrm{d}\P^0\text{-a.e.}
$.
Furthermore, for every $\alpha\in\Ac$ and $\beta\in\Hb(\R^d)$,
\vspace{-2mm}
$$
\hat{h}^{\gm,\hat\alpha^\gm,\beta}_t(Y^\gm_t,Z^\gm_t)
\leq
\hat{h}^{\gm,\hat\alpha^\gm,\hat\beta^\gm}_t(Y^\gm_t,Z^\gm_t)
\leq
\hat{h}^{\gm,\alpha,\hat\beta^\gm}_t(Y^\gm_t,Z^\gm_t)
\quad\mathrm{d}t\times\mathrm{d}\P^0\text{-a.e.}
\vspace{-2mm}
$$
By the comparison theorem for BSDEs (applied to the family $\{\hat{h}^{\gm,\alpha,\beta}\}_{\alpha,\beta}$), we deduce that, for all $t\in[0,T]$,
$
Y^{\gm,\hat\alpha^\gm,\beta}_t
\leq
Y^\gm_t
\leq
Y^{\gm,\alpha,\hat\beta^\gm}_t
\quad\text{a.s.}
$
Hence,
\vspace{-3mm}
$$
Y^\gm_t
=
\essinf_{\alpha\in\Ac}\ \esssup_{\beta\in\Hb(\R^d)} Y^{\gm,\alpha,\beta}_t
=
\esssup_{\beta\in\Hb(\R^d)}\ \essinf_{\alpha\in\Ac} Y^{\gm,\alpha,\beta}_t,
\qquad t\in[0,T].
\vspace{-3mm}
$$
Evaluating at $t=0$ and using the explicit representation of $Y^{\gm,\alpha,\beta}_0$ from the linear BSDE, we obtain $G(\gm)=Y^\gm_0$.

\noindent\textbf{Step 2.}
We now show that $G$ is nondecreasing and can be continuously extended to $\gm=0$ with $G(0)=Y_0$, where $(Y,Z)$ solves ${\rm BSDE}(f,\xi)$.

\noindent Monotonicity in $\gm$ follows directly from the comparison theorem for quadratic BSDEs (Theorem 7.3.1 in \citeauthor{Jianfen_BSDE} \cite{Jianfen_BSDE}), since $\hat{F}^\gm_t(y,z)=f_t(y,z)+\gm|z|^2$ is nondecreasing in $\gm$ for each $(t,y,z)$.

\noindent For $\gm>0$, set
$
U^{0,\gm}:=\frac{Y^\gm-Y}{\gm}\in\S^\infty
\, , \,
V^{0,\gm}:=\frac{Z^\gm-Z}{\gm}\in\Hb.
$
Then $(U^{0,\gm},V^{0,\gm})$ solves ${\rm BSDE}(g^\gm,0)$, where
$
g_t^\gm(u,v)
:=
|Z_t^\gm|^2+l_t^\gm u+m_t^\gm\cdot v,
$
with
\vspace{-2mm}
$$
l^\gm_t
\!\!
:=\frac{f_t(Y_t
\!+
\!
\gm U^{0,\gm}_t
\!,\!
Z_t\!+\!\gm V^{0,\gm}_t)
\!-
\!
f_t(Y_t,\!Z_t \!+\! \gm V^{0,\gm}_t)}{\gm U^{0,\gm}_t}\mathbf{1}_{\{U^{0,\gm}_t\neq0\}}
\, , \,
m^\gm_t:=
\!\!
\int_0^1
\!\!
\!
\partial_z f_t(Y_t,\!Z_t+a\gm V^{0,\gm}_t) \mathrm{d}a.
\vspace{-3mm}
$$
Since $f$ is Lipschitz by Assumption \ref{ass:controlLinf cas controle générales}, there exists $L>0$ such that
\vspace{-2mm}
\begin{equation*}
|l_t^\gm|\le L,
\qquad
|m_t^\gm|\le L,
\qquad
\text{for all } t\in[0,T],\ \gm>0.
\vspace{-2mm}
\end{equation*}
Since $\gm$ ranges in a neighborhood of $0$, standard stability estimates for quadratic BSDEs imply that $(Z^\gm)_\gm$ is bounded in $\H^2$ and $(Y^\gm)_\gm$ is bounded in $\S^2$. Applying $L^p$ estimates to the linear ${\rm BSDE}(g^\gm,0)$, we obtain that for any $p>1$, there exists a constant $C_p>0$ such that, for all $\gm>0$,
\vspace{-6mm}
\begin{equation}\label{ineq:u0etc}
\E\Big[
\sup_{0\le t\le T}|U_t^{0,\gm}|^p
+
\Big(\int_0^T|V_t^{0,\gm}|^2\,\mathrm{d}t\Big)^{p/2}
\Big]
\le C_p.
\vspace{-3mm}
\end{equation}
In particular, $(U^{0,\gm},V^{0,\gm})$ is uniformly bounded in $\S^p\times\H^p$. Together with the stability of quadratic BSDEs as $\gm\to0$, this implies that
as $\gm$ goes to $0$, 
$
\Vert Y^\gm-Y\Vert_{\S^2}
+
\Vert Z^\gm-Z\Vert_{\H^2}
\to0
.$
Hence
$
G(\gm)=Y_0^\gm\to Y_0=:G(0),
$
which gives the continuous extension at $\gm=0$.

\noindent\textbf{Step 3.}
We next show that
$
U_0^{0,\gm}
=
\frac{G(\gm)-G(0)}{\gm}
$
has a finite limit as $\gm\to0$.
Define the limit generator
\vspace{-2mm}
$$
g_t^0(u,v)
:=
|Z_t|^2+\partial_y f_t(Y_t,Z_t)\,u+\partial_z f_t(Y_t,Z_t)\cdot v.
\vspace{-2mm}
$$
Since ${\rm BSDE}(g^0,0)$ is linear with bounded coefficients, it admits a unique solution
$
(U,V)\in\S^\infty\times\Hb(\R^d).
$
Let
$
\delta U^\gm:=U^{0,\gm}-U
$ and $
\delta V^\gm:=V^{0,\gm}-V.
$
Subtracting the BSDEs satisfied by $(U^{0,\gm},V^{0,\gm})$ and $(U,V)$, we see that $(\delta U^\gm,\delta V^\gm)$ solves the linear BSDE
\vspace{-2.5mm}
\begin{equation}\label{eq:d deltaU retravaillé}
\mathrm{d}\delta U^\gm_t
=
\delta V^\gm_t\cdot\mathrm{d}X_t
\!-\!
\Big(
|Z^\gm_t|^2
\!-
\!
|Z_t|^2
+
R^\gm_t
+
\partial_y f_t(Y_t,Z_t)\,\delta U^\gm_t
+
\partial_z f_t(Y_t,Z_t)\cdot\delta V^\gm_t
\Big)\mathrm{d}t,
\vspace{-3mm}
\end{equation}
where
\vspace{-2mm}
$$
R_t^\gm
:=
\frac1\gm
\Big(
f_t(Y_t+\gm U_t^{0,\gm},Z_t+\gm V_t^{0,\gm})
-f_t(Y_t,Z_t)
-\gm U_t^{0,\gm}\,\partial_y f_t(Y_t,Z_t)
-\gm V_t^{0,\gm}\cdot \partial_z f_t(Y_t,Z_t)
\Big).
\vspace{-2.5mm}
$$
By Assumption \ref{ass:controlL2 générales} and a Taylor expansion in $(y,z)$, there exists $C>0$ such that
\vspace{-2mm}
\begin{equation}\label{eq:majo Rh}
|R_t^\gm|
\leq
C\gm\big(|U_t^{0,\gm}|^2+|V_t^{0,\gm}|^2\big).
\vspace{-2mm}
\end{equation}
Together with \eqref{ineq:u0etc}, this implies that $R^\gm\in\H^2$, so that the right-hand side of \eqref{eq:d deltaU retravaillé} is square-integrable.
Since \eqref{eq:d deltaU retravaillé} is a linear BSDE with bounded coefficients, we have
\vspace{-2mm}
\begin{equation*}
\delta U_0^\gm
=
\E\Big[
e^M
\int_0^T e^{N_t}\big(|Z_t^\gm|^2-|Z_t|^2\big)\,\mathrm{d}t
\Big]
+
\E\Big[
e^M
\int_0^T e^{N_t}R_t^\gm\,\mathrm{d}t
\Big],
\vspace{-2mm}
\end{equation*}
where
$
N_t:=\int_0^t \partial_y f_s(Y_s,Z_s)\,\mathrm{d}s
$
and
$
M
:=
{\rm ln} (
\Ec (\int_0^\cdot \partial_z f_s(Y_s,Z_s)\cdot \mathrm{d}X_s)_T
)
$
Since $|\partial_z f_t|\leq L$ by Assumption \ref{ass:controlL2 générales}, we have for all $q\geq1$,
\vspace{-3mm}
$$
\E[e^{qM}]
\leq
\exp\Big(\frac12 Tq^2L^2\Big)
<\infty.
\vspace{-3mm}
$$
Similarly,
$
N_t\leq LT
$
since $|\partial_y f_t|\leq L$.
Using \eqref{eq:majo Rh}, the Cauchy--Schwarz inequality, and the uniform bound \eqref{ineq:u0etc}, we obtain
\vspace{-2mm}
\begin{align*}
|\delta U_0^\gm|
&\leq
e^{LT}\E[e^{2M}]^{1/2}
\Big(
\E\Big[
\Big(\int_0^T \big||Z_t^\gm|^2-|Z_t|^2\big|\,\mathrm{d}t\Big)^2
\Big]^{1/2}
\\[-1mm]
&\hspace{2.2cm}
+
\sqrt{2}C\gm
\E\Big[
\Big(\int_0^T |U_t^{0,\gm}|^2\,\mathrm{d}t\Big)^2
+
\Big(\int_0^T |V_t^{0,\gm}|^2\,\mathrm{d}t\Big)^2
\Big]^{1/2}
\Big).
\vspace{-2mm}
\end{align*}
\vspace{-5mm}

\noindent Since $Z^\gm\to Z$ in $\H^{2p}$, the first term converges to $0$ as $\gm\to0$, while the second term is of order $O(\gm)$ by \eqref{ineq:u0etc}. Hence
$
|\delta U_0^\gm|\to0
$
as $\gm\to0$, and therefore
\vspace{-2mm}
$$
\frac{G(\gm)-G(0)}{\gm}
=
U_0^{0,\gm}
=
U_0+\delta U_0^\gm
\longrightarrow U_0
\qquad \text{as }\gm\to0.
\vspace{-2mm}
$$
Thus $G$ is differentiable at $0$, with
$
G'(0)=U_0.
$
Finally, solving the linear ${\rm BSDE}(g^0,0)$ explicitly yields
$
U_0
=
\E^{\lm^*}
[
\int_0^T \Kc_t^* |Z_t|^2\,\mathrm{d}t
],
$
where $(\lm^*,k^*,\Kc^*)$ are defined in \eqref{eqdef:kopt_and_lambdaopt}. In other words,
$
G'(0)=\Vert Z\Vert_{\L^2(\P^{\lm^*},k^*)}^2.
$
This completes the proof.
\ep

\noindent We now establish that, for $r \geq 0$, $\tl V_2(r)=\overline{V}_2(r)$, where $\tl V_2$ is defined in \eqref{eqdef:tl V(r)}. 
Fix $\alpha \in \Ac$, and let $\hat{J}_\gm$ be defined by \eqref{eqdef:hatJ}.
Our objective is to use Lemma \ref{lem:strong duality} to prove the min--max identity
\vspace{-3.5mm}
$$
\inf_{\gm>0} G^\alpha(\gm)
=
\sup_{\beta\in\Hb(\R^d)} \inf_{\gm>0} \hat{J}_\gm(\alpha,\beta),
\,\, \text{
where}
\,\, 
G^\alpha(\gm)
:=
\sup_{\beta\in\Hb(\R^d)} \hat{J}_\gm(\alpha,\beta).
\vspace{-2.5mm}
$$

\begin{Remark}
{\rm
The choice of the deviation set $\Bc^2_\alpha$ is crucial for our results and computations.
When we use Lagrange multipliers to dualize the constraint, the resulting Lagrangian can still be written as a single expectation under $\P^{\lm^\alpha+\beta}$. 
With a different criterion, such a representation would fail, and we would lose the dynamic programming structure. 
In particular, the dual problem $G^\alpha$ could no longer be represented as the initial value of a BSDE.
}
\end{Remark}
\vspace{-4mm}

\begin{Proposition}\label{prop:diff of Galpha(gm)}
Let Assumption \ref{ass:controlL2 générales} hold, and let $\xi\in\L^\infty$.
Then $G^\alpha$ is a nondecreasing continuous function, which admits a continuous extension at $\gm=0$.
Moreover, $G^\alpha$ is differentiable, with
\vspace{-3mm}
$$
(G^\alpha)'(\gm)
=
\E^{\lm^\alpha+2\gm Z^{\gm,\alpha}}
\Big[
\int_0^T \Kc_t^\alpha |Z^{\gm,\alpha}_t|^2\,\mathrm{d}t
\Big],
\vspace{-2mm}
$$
where $(Y^{\gm,\alpha},Z^{\gm,\alpha})$ denotes the solution of ${\rm BSDE}(\hat{f}^{\gm,\alpha},\xi)$.
\end{Proposition}

\vspace{-4mm}

\begin{Remark}
{\rm
Although Propositions \ref{prop:diff of Galpha(gm)} and \ref{prop:diff tl V2} have similar statements, their proofs rely on different structural properties of the underlying BSDEs.
For Proposition \ref{prop:diff of Galpha(gm)}, we exploit the linear structure (for fixed $\alpha$) to obtain differentiability of $G^\alpha$ for all $\gm>0$.
In contrast, in the proof of Proposition \ref{prop:diff tl V2}, we only obtain differentiability at $\gm=0$, since the quadratic term vanishes at this point.
}
\end{Remark}

\vspace{-3.5mm}

\proof 
\textbf{Step $1$.}
We first provide a BSDE representation for $G^\alpha$. 
Consider the generator
\vspace{-2mm}
$$
\hat{f}^{\gm,\alpha}_t(y,z)
:=
\sup_{b\in\R^d}
\big\{
f^\alpha_t(y,z)-b\cdot z-\tfrac{1}{4\gm}|b|^2
\big\}
=
f^\alpha_t(y,z)+\gm|z|^2,
\vspace{-3mm}
$$
where we recall that $f^\alpha_t(y,z):=l^\alpha_t-k^\alpha_ty-\lm^\alpha_t\cdot z$.
For all $\gm\ge0$, the quadratic ${\rm BSDE}(\hat{f}^{\gm,\alpha},\xi)$ is well posed and admits a unique solution $(Y^{\gm,\alpha},Z^{\gm,\alpha})\in\S^\infty\times\Hb(\R^d)$, by Assumption \ref{ass:controlL2 générales} and the well-posedness theorem for quadratic BSDEs (Theorem 7.3.3 in \citeauthor{Jianfen_BSDE} \cite{Jianfen_BSDE}), together with the estimate in Theorem 7.2.1.
Moreover, by adapting the first step of the proof of Proposition \ref{prop:diff tl V2}, and introducing
$
\hat\beta^\gm:=-2\gm Z^{\gm,\alpha}\in\Hb(\R^d),
$
one obtains that, for all $t\in[0,T]$,
$
Y^{\gm,\alpha}_t
=
Y^{\gm,\alpha,\hat\beta^\gm}_t
\quad\text{a.s.}
$ and $
Z^{\gm,\alpha}_t
=
Z^{\gm,\alpha,\hat\beta^\gm}_t
\quad \mathrm{d}t\times\mathrm{d}\P^0\text{-a.e.}
.$
Furthermore, 
\vspace{-4mm}
$$
Y^{\gm,\alpha}_0
=
Y^{\gm,\alpha,\hat\beta^\gm}_0
=
\sup_{\beta\in\Hb(\R^d)} Y^{\gm,\alpha,\beta}_0
=
G^\alpha(\gm).
\vspace{-2mm}
$$
Adapting once again the arguments of Proposition \ref{prop:diff tl V2} shows that $G^\alpha$ is nondecreasing and continuous, and admits a continuous extension at $\gm=0$ given by $G^\alpha(0):=Y_0^\alpha$, where $(Y^\alpha,Z^\alpha)$ solves ${\rm BSDE}(f^\alpha,\xi)$.

\noindent\textbf{Step $2$.} We prove the differentiability of $\gm\mapsto \tl Y^{\gm,\alpha}_0$.
Define
$
\tl Y^{\gm,\alpha}
:=
\exp(2\gm Y^{\gm,\alpha})
$ and $
\tl Z^{\gm,\alpha}
:=
2\gm \tl Y^{\gm,\alpha} Z^{\gm,\alpha}.
$
Fix $\gm>0$. By Itô's formula, $(\tl Y^{\gm,\alpha},\tl Z^{\gm,\alpha})$ solves ${\rm BSDE}(\tl f,\tl\xi_\gm)$ with
\vspace{-1mm}
$$
\tl f_t(y,z)
:=
2\gm l^\alpha_t y
-
k^\alpha_t y\ln(y)
-
\lm^\alpha_t\cdot z,
\qquad
\tl\xi_\gm:=\exp(2\gm\xi).
\vspace{-1mm}
$$
From Theorem 7.2.1 in \citeauthor{Jianfen_BSDE} \cite{Jianfen_BSDE}, for all $|h|<1\wedge\gm/2$, we obtain the uniform bound
$
\Vert Y^{\gm+h,\alpha}\Vert_\infty\le M^\gm,
$
for some constant $M^\gm>0$. It follows that there exist constants $0<c^\gm\le C^\gm<\infty$ such that
\vspace{-1mm}
$$
c^\gm\le\tl Y_t^{\gm+h,\alpha}\le C^\gm
\quad \text{for all }t\in[0,T],\ \textit{a.s.}
\vspace{-1mm}
$$
Let $\varphi:\R\to\R$ be a $C^1$ and  $K$-Lipschitz function such that
$
\varphi(x)=x {\rm ln}(x)
$
for all $x\in[c^\gm/2,2C^\gm]$.
Then $(\tl Y^{\gm+h,\alpha},\tl Z^{\gm+h,\alpha})$ is the unique solution of ${\rm BSDE}(H^{\gm+h,\varphi},\tl\xi_{\gm+h})$, where
\vspace{-1.5mm}
$$
H_t^{\gm+h,\varphi}(y,z)
:=
2(\gm+h)l_t^\alpha y
-
k_t^\alpha\varphi(y)
-
\lm_t^\alpha\cdot z,
\qquad
\tl\xi_{\gm+h}:=\exp(2(\gm+h)\xi).
\vspace{-1.5mm}
$$
In particular, ${\rm BSDE}(H^{\gm+h,\varphi},\tl\xi_{\gm+h})$ is well posed, with unique solution $(\tl Y^{\gm+h,\alpha},\tl Z^{\gm+h,\alpha})\in\S^2\times\H^2(\R^d)$.
Define the difference quotients
$
\tl U^{\gm,\alpha,h}
:=
\frac{\tl Y^{\gm+h,\alpha}-\tl Y^{\gm,\alpha}}{h},
$
and
$
\tl V^{\gm,\alpha,h}
:=
\frac{\tl Z^{\gm+h,\alpha}-\tl Z^{\gm,\alpha}}{h}.
$
Then $(\tl U^{\gm,\alpha,h},\tl V^{\gm,\alpha,h})$ solves ${\rm BSDE}(\tl g^{\gm,\alpha,h},\zeta^{\gm,h})$, with generator
\vspace{-1mm}
\begin{equation*}
\tl g_t^{\gm,\alpha,h}(u,v)
:=
\frac{
2(\gm+h)l_t^\alpha(\tl Y_t^{\gm,\alpha}+hu)
-
2\gm l_t^\alpha\tl Y_t^{\gm,\alpha}
}{h}
-
\lm_t^\alpha\cdot v
-
k_t^\alpha
\frac{
\varphi(\tl Y_t^{\gm,\alpha}+hu)-\varphi(\tl Y_t^{\gm,\alpha})
}{h},
\vspace{-2mm}
\end{equation*}
and terminal condition
$
\zeta^{\gm,h}
:=
\frac{\exp(2(\gm+h)\xi)-\exp(2\gm\xi)}{h}.
$
As $h\to0$, we have
\vspace{-2.5mm}
$$
\tl g_t^{\gm,\alpha,h}(u,v)
\xrightarrow[h\to0]{\textit{a.s.}}
\tl g_t^{\gm,\alpha}(u,v)
:=
2l_t^\alpha \tl Y_t^{\gm,\alpha}
+
2\gm l_t^\alpha u
-
\lm_t^\alpha\cdot v
-
k_t^\alpha\varphi'(\tl Y_t^{\gm,\alpha})u,
\vspace{-2mm}
$$
and, since $\xi\in\L^\infty$,
$
\E[|\zeta^{\gm,h}-\zeta^\gm|^2]\xrightarrow[h\to0]{}0,
$
where
$
\zeta^\gm:=2\xi\exp(2\gm\xi).
$
By the stability theorem for Lipschitz BSDEs (Theorem 4.4.3 in \citeauthor{Jianfen_BSDE} \cite{Jianfen_BSDE}), we obtain
\vspace{-2.5mm}
\begin{equation*}
\E\Big[
\sup_{0\le t\le T}|\tl U_t^{\gm,\alpha}-\tl U_t^{\gm,\alpha,h}|^2
+
\int_0^T|\tl V_t^{\gm,\alpha}-\tl V_t^{\gm,\alpha,h}|^2\,\mathrm{d}t
\Big]
\xrightarrow[h\to0]{}0,
\vspace{-3mm}
\end{equation*}
where $(\tl U^{\gm,\alpha},\tl V^{\gm,\alpha})$ is the unique solution of ${\rm BSDE}(\tl g^{\gm,\alpha},\zeta^\gm)$.

\noindent\textbf{Step $3$.}
We now prove the differentiability of $G^\alpha$. This follows from the differentiability of $\gm\mapsto \tl Y_0^{\gm,\alpha}$, since
$
G^\alpha(\gm)=Y_0^{\gm,\alpha}={\rm ln}(\tl Y_0^{\gm,\alpha})/(2\gm).
$
By the chain rule,
$
(G^\alpha)'(\gm)=U_0^{\gm,\alpha},
$
where
$
U^{\gm,\alpha}
:=
-\frac{1}{2\gm^2}{\rm ln}(\tl Y^{\gm,\alpha})
+
\frac{\tl U^{\gm,\alpha}}{2\gm\tl Y^{\gm,\alpha}}.
$
Applying Itô's formula to $U^{\gm,\alpha}$ and rearranging terms, we obtain
\vspace{-2.5mm}
\begin{equation*}
\mathrm{d}U_t^{\gm,\alpha}
=
V_t^{\gm,\alpha}\cdot\mathrm{d}X_t
-
\Big(
|Z_t^{\gm,\alpha}|^2
-
k_t^\alpha U_t^{\gm,\alpha}
-
(\lm_t^\alpha-2\gm Z_t^{\gm,\alpha})\cdot V_t^{\gm,\alpha}
\Big)\mathrm{d}t,
\vspace{-2mm}
\end{equation*}
where
$
V^{\gm,\alpha}
:=
-\frac{Z^{\gm,\alpha}}{\gm}
-
\frac{\tl U^{\gm,\alpha} Z^{\gm,\alpha}}{\tl Y^{\gm,\alpha}}
+
\frac{\tl V^{\gm,\alpha}}{2\gm\tl Y^{\gm,\alpha}}.
$
Thus, $(U^{\gm,\alpha},V^{\gm,\alpha})$ solves a linear BSDE with bounded coefficients under the measure $\P^{\lm^\alpha+2\gm Z^{\gm,\alpha}}$, which yields the desired representation
\vspace{-2mm}
$$
(G^\alpha)'(\gm)
=
U_0^{\gm,\alpha}
=
\E^{\lm^\alpha+2\gm Z^{\gm,\alpha}}
\Big[
\int_0^T \Kc_t^\alpha |Z_t^{\gm,\alpha}|^2\,\mathrm{d}t
\Big].
$$
\ep

\begin{Proposition}\label{prop:saddle-point Lagrangien alpha fixé}
Let Assumption \ref{ass:controlL2 générales} hold, and let $\xi \in \L^\infty$. Then, for all $r>0$, we have
\vspace{-3mm}
$$
\displaystyle\inf_{\gm>0}\displaystyle\sup_{\beta\in\Hb(\R^d)}\Lc_r^\alpha(\gm,\beta)
=
\displaystyle\sup_{\beta\in\Hb(\R^d)}\displaystyle\inf_{\gm>0}\Lc_r^\alpha(\gm,\beta).
\vspace{-3mm}
$$
\end{Proposition}

\proof 
Let $\Xc := \Hb(\R^d)$, $
\psi(\beta)
:=
\E^{\lm^\alpha+\beta}
[\Kc_T^\alpha \xi+\int_0^T \Kc_t^\alpha l^\alpha_t\,\mathrm{d} t],
$ and set 
\vspace{-2mm}
$$
\varphi(\beta)
:=
\E^{\lm^\alpha+\beta}
\Big[\int_0^T \Kc_t^\alpha |\beta_t|^2 \mathrm{d} t\Big]
\,\, , \,\, 
H(\lm)
:=
\sup_{ \beta \in \Xc} \big\{\psi(\beta) - \lm \varphi(\beta)\big\}
=
G^{\alpha}\Big(\frac{1}{4\lm}\Big).
\vspace{-2mm}
$$
By Proposition \ref{prop:diff of Galpha(gm)} and the BSDE representation of $G^{\alpha}$, it is clear that $H$ is bounded. Indeed, since $\xi \in \L^{\infty}$ and $k,l$ are bounded, we have
\vspace{-1mm}
$$
\vert \xi \vert_{\infty} e^{T \vert k^\alpha \vert_{\infty}} + 
T \vert l^\alpha \vert_{\infty} e^{T \vert k^\alpha \vert_{\infty}}
\geq 
G^{\alpha}(\gm)
\geq 
\psi(0) - \lm \varphi(0) = \psi(0),
\vspace{-1mm}
$$
which is finite by Assumption \ref{ass:controlL2 générales}.
By Step 1 in the proof of Proposition \ref{prop:diff of Galpha(gm)}, for all $\lm>0$, there exists an element $\beta_\lm \in \Xc$ such that
\vspace{-2mm}
$$
H(\lm)
=
\sup_{ \beta \in \Xc} \big\{\psi(\beta) - \lm \varphi(\beta)\big\}
=
\psi(\beta_\lm) - \lm \varphi(\beta_\lm),
\vspace{-2mm}
$$
and, furthermore, $\lm \mapsto H(\lm)$ is differentiable with
$
H'(\lm) = -\varphi(\beta_\lm).
$
Hence the assumptions of Lemma \ref{lem:strong duality} are satisfied, and we may apply it. 
\ep 

\smallskip 

\noindent \textit{Proof of Theorem \ref{thm:V'(0) Linf et L2} \ref{thm:V'(0) L2}.} We first prove that $\overline{V}_2(r)=\tl V_2(r)$. For $r=0$, since $G$ is nondecreasing, we have 
$
\tl V_2(0)=\displaystyle\inf_{\gm>0}G(\gm)=G(0)=Y_0=\overline{V}_2(0).
$
For $r>0$, Proposition \ref{prop:saddle-point Lagrangien alpha fixé} yields 
\vspace{-3mm}
\begin{equation*}
\displaystyle\inf_{\alpha\in\Ac}\displaystyle\inf_{\gm>0}\displaystyle\sup_{\beta\in\Hb(\R^d)}\Lc_r^\alpha(\gm,\beta)
=
\displaystyle\inf_{\alpha\in\Ac}\displaystyle\sup_{\beta\in\Hb(\R^d)}\displaystyle\inf_{\gm>0}\Lc_r^\alpha(\gm,\beta).
\end{equation*}
By inverting $\inf_{\alpha\in\Ac}$ and $\inf_{\gm>0}$, we get $\tl V_2(r)= \overline{V}_2(r)$.
Thus, we have	
$
\overline{V}_2(r)
=
\inf_{\gm>0}\{G(\gm)+\frac{r^2}{4\gm}\},	
$
where $G:\R_+\to\R$ is nondecreasing and differentiable at $\gm=0$. 
Applying Lemma \ref{lem:concluding lemma}, we obtain $\overline{V}_2'(0)$ as in Theorem \ref{thm:V'(0) Linf et L2} \ref{thm:V'(0) L2}, which completes the proof. \ep 

\medskip 

\begin{Remark}\label{rem:foealphaL2}
We see here that a first-order expansion of
$
\alpha_r^\ast
=
\argmin_{\alpha \in \Ac}
\sup_{\beta\in\Bc^2_\alpha(r)}
J(\alpha,\beta)
$
is difficult to obtain. Indeed, the distributionally robust feedback is
$
\alpha_r^\ast
=
\alpha^\ast\!(Y^{\gm_r},Z^{\gm_r}),
$
where
$
\gm_r
=
\argmin_{\gm>0}
\{
G(\gm)+\frac{r^2}{4\gm}
\},
$
and $(Y^{\gm},Z^{\gm})$ denotes the solution to the quadratic ${\rm BSDE}(\hat{F}^\gm,\xi)$.
However, the existence of
\(
\gm_r
=
\argmin_{\gm>0} \{G(\gm)+\frac{r^2}{4\gm} \}
\)
is not guaranteed, and we would additionally need a first-order expansion at $r=0$ of this family.
This would require extra regularity of $G$, which we were not able to establish.
\end{Remark}

\section{Sensitivities for Mixed Optimal Control and Stopping}\label{proof thm mixed}

\subsection{Deriving the $\L^\infty$ Sensitivity of DROCS Problems}
We begin with the $\L^\infty$ sensitivity.
For $r\geq0$ and $\beta\in\Bc^\infty(r)$, let the generator $f, f^\alpha, 
h^{\alpha, \beta}$, and $F^r$ be defined by \eqref{eqdef:genLinftyopticont}.

\noindent \textit{Proof of Theorem \ref{thm:sensimixedproblem} \ref{thm:Vup=Vdown mixed inf}.}
\textbf{Step $1$.} We prove that 
\vspace{-2.5mm}
\begin{equation}\label{eq:V(r)=Yr} 
\overline{\V}_\infty(r)= \underline{\V}_\infty(r)=Y^r_0.
\vspace{-2.5mm}
\end{equation}
Let $( S^{\alpha, \beta} (\tau, \xi), \pi^{\alpha, \beta}(\tau, \xi) ) $ be the solution to ${\rm BSDE}_{\tau} ( h^{\alpha, \beta} , \xi) $ (as defined in Equation \eqref{eqdef:BSDE random term time}).
We first establish that
 \vspace{-3mm}
 \begin{equation*}
    \inf_{\alpha \in \Ac} \inf_{\tau\in\Tc} 
    \sup_{\beta\in\Bc^\infty(r)} S^{\alpha, \beta}_0(\tau,\xi) = Y^r_0
    =
    \sup_{\beta\in\Bc^\infty(r)}\inf_{\alpha \in \Ac} \inf_{\tau\in\Tc} S^{\alpha,\beta}_0(\tau,\xi),
    \vspace{-3mm}
  \end{equation*}
  where $ (Y^r, Z^r, K^r)$ is the unique solution of ${\rm RBSDE}(F^r, \xi)$.
Fix $\beta\in\Bc^\infty(r)$ and $\tau \in\Tc$. 
Since $\tau$ is a bounded stopping time and $h^{\alpha, \beta}$ is Lipschitz, by Proposition 4.1.2 and Remark 4.3.2 in \citeauthor{Jianfen_BSDE} \cite{Jianfen_BSDE}, ${\rm BSDE}_{\tau}(h^{\alpha, \beta}, \xi)$ is well-posed, and we have 
$
\E^{{\lm^\alpha + \beta}}[\Kc_\tau^\alpha \xi_{\tau} +\int_0^\tau \Kc_t^\alpha l^\alpha_t\mathrm{d}t] 
 = S_0^{\alpha, \beta}(\tau,\xi)$. 
Likewise, by Theorem 6.3.1 of \citeauthor{Jianfen_BSDE} \cite{Jianfen_BSDE}, ${\rm RBSDE}(F^r,\xi)$ is well-posed.
We verify that the family $(F^r,h^{\alpha,\beta})_{\alpha\in\Ac,\ \beta\in\Bc^\infty(r)}$, defined in \eqref{eqdef:genLinftyopticont}, forms a family of standard drivers in the sense of \citeauthor{el1997reflected} \cite{el1997reflected}.
Therefore, Theorem 5.1 of \cite{el1997reflected} applies and yields the desired equality \eqref{eq:V(r)=Yr}.

\noindent \textbf{Step $2$.} Let $( Y, Z, K) \in\S^2\times\H^2(\R^d) \times \mathbb{A}^2$ be the unique solution of ${\rm RBSDE}(f , \xi)$, which exists under the conditions of Theorem \ref{thm:sensimixedproblem}. 
We now show that, as $r\to0$, one has $\tau^r\to \tl\tau$ \textit{a.s.}, where $\tilde{\tau}$ is defined in \eqref{eqdef:tau_opt_mixed} and
$
  \tau^r:=\inf\{t\in[0,T]:Y^r_t=\xi_t\}.
$
First, for $0\leq r\leq r'$, we have $
    F^r_t(y, z)\leq F^{r'}_t(y, z)$, which, by the comparison theorem for RBSDEs (Theorem 6.2.5 in \citeauthor{Jianfen_BSDE} \cite{Jianfen_BSDE}), implies that for all $t\in[0,T]$, $ 
    Y^{r}_t\leq Y^{r'}_t\leq \xi_t \quad \textit{a.s.}
$
Furthermore, since $Y^{r'} \leq \xi$, for $t=\tau^r$, $Y^{r'}_{\tau^r} \leq \xi_{\tau^r}$ and therefore 
 \vspace{-3mm}
  \begin{equation}\label{eq:tau r croissant}
    \tau^{r'}\leq\tau^r \quad \textit{a.s.} 
     \vspace{-3.5mm}
  \end{equation}
This implies the existence of $\tau^*\leq\tl{\tau}$ such that \textit{a.s.} $
    \tau^r\rightarrow \tau^*$. 
By Assumption \ref{ass:controlLinf cas controle générales}, $f$ is Lipschitz. 
Therefore, by Theorem 6.2.3 in \citeauthor{Jianfen_BSDE} \cite{Jianfen_BSDE}, there exists $M>0$ such that 
\vspace{-2.5mm}
  \begin{equation}\label{ineq:estimRBSDELinfyt}
\Vert Y^r - Y \Vert_{\S^2}^2 + \Vert K^r - K \Vert_{\S^2}^2 + \Vert Z^r - Z \Vert_{\H^2}^2
\leq M r^2.
\vspace{-2mm}
  \end{equation}
It follows that $Y^r_{\tau^r} - Y_{\tau^r} \to 0$ in probability. 
Using the continuity of $\xi$ on $[0,T)$, together with $\xi_{T^-}\geq \xi_T$, one deduces that \textit{a.s.} $
    \xi_{\tau^r}\rightarrow\xi_{\tl\tau}\mathbf{1}_{\{\tl\tau<T\}}+\xi_{T^-}\mathbf{1}_{\{\tl\tau=T\}}$. 
Furthermore, as a solution to an RBSDE, $Y$ is continuous, and therefore \textit{a.s.}, $Y_{\tau^r}\rightarrow Y_{\tl\tau}$.
Since $Y$ is continuous and satisfies $Y_T=\xi_T$, it follows that $Y_{\tau^*}=\xi_{\tau^*}$, hence $\tau^*\geq \tl{\tau}$. 
Combining this with \eqref{eq:tau r croissant}, we obtain the claim.

\noindent \textbf{Step $3$.}
We next prove that the limit
$
\lim_{r\to0}\frac{Y_0^r-Y_0}{r}
$
exists, and we identify its value. Let $(U,V)$ denote the unique solution of the linear ${\rm BSDE}_{\tl{\tau}}(g^0,0)$, where
\vspace{-2mm}
$$
g_t^0(u,v)
:=
\partial_y f(Y_t,Z_t)u+\partial_z f(Y_t,Z_t)\cdot v+|Z_t|.
\vspace{-2mm}
$$
For notational simplicity, we omit the dependence on $r$, define
$
\delta Y_t^r
:=
\frac{Y_{t\wedge\tl{\tau}}^r-Y_{t\wedge\tl{\tau}}}{r}
$, $
\delta Z_t^r
:=
\frac{Z_{t\wedge\tl{\tau}}^r-Z_{t\wedge\tl{\tau}}}{r}
$
,
$
\dU_t:=\delta Y_t^r-U_{t\wedge\tl{\tau}}
$ and $
\dV_t:=\delta Z_t^r-V_{t\wedge\tl{\tau}}.
$
For $r>0$, let
\vspace{-2mm}
$$
g_t^r(u,v)
:=
\frac{f_t(Y_t+ru,Z_t+rv)-f_t(Y_t,Z_t)}{r}
+
|Z_t+rv|.
\vspace{-2mm}
$$
Since $K$ is continuous and satisfies the Skorokhod condition, one has $K_t=0$ for all $t\le \tl{\tau}$.
Therefore, for $t\in[0,\tl{\tau}]$, letting 
$
\Delta_t^r
:=
g_t^r(\delta Y_t^r,\delta Z_t^r)-g_t^0(U_t,V_t)
$, we have 
\vspace{-2.5mm}
$$
\mathrm{d}\dU_t
=
\dV_t\cdot \mathrm{d}X_t
-
\Delta_t^r\,\mathrm{d}t
+
\frac1r\,\mathrm{d}K_t^r.
\vspace{-2.5mm}
$$ 
\noindent Using Itô's formula, and since by Remark \ref{rem:S2*H2=mtg}
$
\int_{0}^{t}\dU_s\dV_s\cdot \mathrm{d}X_s
$ is a martingale, we get
\vspace{-2.5mm}
\begin{equation}\label{eq:ito deltaU2}
\E\Big[
(\dU_t)^2+\int_{t\wedge\tl{\tau}}^{\tl{\tau}}|\dV_s|^2\mathrm{d}s
\Big]
=
I_t^{1,r}+I_t^{2,r}
\,\, \text{ for $0\le t\le T$,}
\vspace{-2.5mm}
\end{equation}
$
I_t^{1,r}
:=
\E [
-\frac2r\int_{t\wedge\tl{\tau}}^{\tl{\tau}}\dU_s\,\mathrm{d}K_s^r
]
$, $
I_t^{2,r}
:=
\E [
-2\int_{t\wedge\tl{\tau}}^{\tl{\tau}}\dU_s\Delta_s^r\,\mathrm{d}s
]
$
.We now estimate separately the terms $I_t^{1,r}$ and $I_t^{2,r}$.

\smallskip

\noindent \textit{Estimate of $I_t^{1,r}$.}
Since $K^r$ is continuous and nondecreasing, and since $Y_t^r<\xi_t$ on $[0,\tau^r)$, we have that for $ 0\le t\le \tau^r$, 
$
K_t^r=0 \textit{a.s.}
$
Moreover, since $Y\le \xi$, we have
$
\int_t^{\tl{\tau}}(Y_s-\xi_s)\,\mathrm{d}K_s^r\le 0.
$
Therefore, for $t\in[0,\tl{\tau}]$,
\vspace{-2mm}
\begin{equation}\label{eqdef:estimDKint}
\begin{split}
-\int_t^{\tl{\tau}}\dU_s\,\mathrm{d}K_s^r
&=
\int_t^{\tl{\tau}}U_s\,\mathrm{d}K_s^r
-\frac1r\int_t^{\tl{\tau}}Y_s^r\,\mathrm{d}K_s^r
+\frac1r\int_t^{\tl{\tau}}Y_s\,\mathrm{d}K_s^r
\\
&=
\int_t^{\tl{\tau}}U_s\,\mathrm{d}K_s^r
+\frac1r\int_t^{\tl{\tau}}(Y_s-\xi_s)\,\mathrm{d}K_s^r
\\
&\le
\int_t^{\tl{\tau}}U_s\,\mathrm{d}K_s^r
=
\int_{\tau^r}^{\tl{\tau}}U_s\,\mathrm{d}K_s^r.
\end{split}
\end{equation}
\vspace{-3mm}

\noindent 
Now, integrating by parts and using $U_{\tl{\tau}}=0$ together with $K_{\tau^r}^r=0$, we get
\vspace{-2mm}
$$
\int_{\tau^r}^{\tl{\tau}}U_s\,\mathrm{d}K_s^r
=
-\int_{\tau^r}^{\tl{\tau}}K_s^r\,\mathrm{d}U_s
=
-\int_{\tau^r}^{\tl{\tau}}K_s^rV_s\cdot \mathrm{d}X_s
+
\int_{\tau^r}^{\tl{\tau}}K_s^r g_s^0(U_s,V_s)\,\mathrm{d}s.
\vspace{-2mm}
$$
Taking expectations in \eqref{eqdef:estimDKint} and using that the stochastic integral is a true martingale, we obtain
\vspace{-2mm}
$$
I_t^{1,r}
\le
\frac2r
\E\Big[
\int_{\tau^r}^{\tl{\tau}}|K_s^r|\,|g_s^0(U_s,V_s)|\,\mathrm{d}s
\Big].
\vspace{-2mm}
$$
Applying Cauchy--Schwarz yields
\vspace{-2mm}
\begin{align*}
I_t^{1,r}
&\le
\frac2r
\E\Big[\int_{\tau^r}^{\tl{\tau}}|K_s^r|^2\,\mathrm{d}s\Big]^{1/2}
\E\Big[\int_{\tau^r}^{\tl{\tau}}|g_s^0(U_s,V_s)|^2\,\mathrm{d}s\Big]^{1/2}
\\
&\le
\frac2r
\E\Big[(\tl{\tau}-\tau^r)\sup_{\tau^r\le s\le \tl{\tau}}|K_s^r|^2\Big]^{1/2}
\E\Big[\int_{\tau^r}^{\tl{\tau}}|g_s^0(U_s,V_s)|^2\,\mathrm{d}s\Big]^{1/2}.
\end{align*}
\vspace{-4mm}

\noindent 
Since $K_s=0$ for $s\le \tl{\tau}$, we have on $[\tau^r,\tl{\tau}]$,
$
K_s^r=K_s^r-K_s.
$
Hence, using \eqref{ineq:estimRBSDELinfyt},
\vspace{-2mm}
\begin{align*}
I_t^{1,r}
&\le
\frac{2\sqrt{T}}r
\E\Big[\sup_{\tau^r\le s\le \tl{\tau}}|K_s^r-K_s|^2\Big]^{1/2}
\E\Big[\int_{\tau^r}^{\tl{\tau}}|g_s^0(U_s,V_s)|^2\,\mathrm{d}s\Big]^{1/2}
\\
&\le
\frac{2\sqrt{T}}r
\E\Big[\sup_{0\le s\le T}|K_s^r-K_s|^2\Big]^{1/2}
\E\Big[\int_{\tau^r}^{\tl{\tau}}|g_s^0(U_s,V_s)|^2\,\mathrm{d}s\Big]^{1/2}
\\
&\le
2\sqrt{TC}\,
\E\Big[\int_{\tau^r}^{\tl{\tau}}|g_s^0(U_s,V_s)|^2\,\mathrm{d}s\Big]^{1/2}.
\end{align*}
\vspace{-4mm}

\noindent Now, since $(U,V)\in\S^2\times\H^2$, Assumption \ref{ass:controlLinf cas controle générales} implies that
$
\int_0^T|g_s^0(U_s,V_s)|^2\,\mathrm{d}s\in\L^1.
$
Moreover, since $\tau^r\xrightarrow[r\to0]{\textit{a.s.}}\tl{\tau}$, we have
$
\int_{\tau^r}^{\tl{\tau}}|g_s^0(U_s,V_s)|^2\,\mathrm{d}s
\xrightarrow[r\to0]{\textit{a.s.}}
0
$.
By dominated convergence,
$
\E\Big[\int_{\tau^r}^{\tl{\tau}}|g_s^0(U_s,V_s)|^2\,\mathrm{d}s\Big]
\xrightarrow[r\to0]{}0,
$
and therefore
$
\sup_{0\le t\le T}I_t^{1,r}\xrightarrow[r\to0]{}0.
$

\smallskip

\noindent \textit{Estimate of $I_t^{2,r}$.}
Since $f$ is Lipschitz by Assumption \ref{ass:controlLinf cas controle générales}, there exists $C>0$ such that
$
|\Delta_s^r|
\le
C|\dU_s|+|g_s^r(U_s,V_s)-g_s^0(U_s,V_s)|.
$
Hence, by Young's inequality 
\vspace{-0mm}
\begin{align*}
I_t^{2,r}
&\le
\E\Big[
\int_{t\wedge\tl{\tau}}^{\tl{\tau}}2|\dU_s|\,|\Delta_s^r|\,\mathrm{d}s
\Big]
\\
&\le
(2C+1)\E\Big[
\int_{t\wedge\tl{\tau}}^{\tl{\tau}}|\dU_s|^2\,\mathrm{d}s
\Big]
+
\E\Big[
\int_0^{\tl{\tau}}|g_s^r(U_s,V_s)-g_s^0(U_s,V_s)|^2\,\mathrm{d}s
\Big].
\end{align*}
\vspace{-5mm}

\noindent Combining the above estimates in \eqref{eq:ito deltaU2}, we obtain for all $0\le t\le T$,
\vspace{-3mm}
\begin{align*}
\E\Big[
(\dU_t)^2+\int_{t\wedge\tl{\tau}}^{\tl{\tau}}|\dV_s|^2\,\mathrm{d}s
\Big]
&\le
(2C+1)
\E\Big[
\int_{t\wedge\tl{\tau}}^{\tl{\tau}}|\dU_s|^2\,\mathrm{d}s
\Big]
\\
&\quad
+
\E\Big[
\int_0^{\tl{\tau}}|g_s^r(U_s,V_s)-g_s^0(U_s,V_s)|^2\,\mathrm{d}s
\Big]
\\
&\quad
+
2\sqrt{TC}\,
\E\Big[
\int_{\tau^r}^{\tl{\tau}}|g_s^0(U_s,V_s)|^2\,\mathrm{d}s
\Big]^{1/2}.
\end{align*}
\vspace{-5mm}

\noindent 
Finally, 
$
g_s^r(U_s,V_s)\rightarrow g_s^0(U_s,V_s)
$ $\mathrm{d}s\otimes\mathrm{d}\P^0\textit{-a.e.}$ 
and
$
|g_s^r(U_s,V_s)-g_s^0(U_s,V_s)|
\le
C(|U_s|+|V_s|),
$
and another application of the dominated convergence theorem yields
$
\E [
\int_0^{\tl{\tau}}|g_s^r(U_s,V_s)-g_s^0(U_s,V_s)|^2\,\mathrm{d}s
]
\rightarrow0.
$
Applying the Gronwall inequality, we conclude that
\vspace{-2mm}
$$
\E\Big[
(\dU_t)^2+\int_{t\wedge\tl{\tau}}^{\tl{\tau}}|\dV_s|^2\,\mathrm{d}s
\Big]
\xrightarrow[r\to0]{}0.
\vspace{-1.5mm}
$$
This concludes the proof, since
$
U_0
=
\E^{\tl{\lm}^*}
[
\int_0^{\tl{\tau}}\Kc_s^*|Z_s|\,\mathrm{d}s
].
$
\ep

\begin{Remark}\label{rem:S2*H2=mtg}
{\rm 
  If $U\in\S^2$ and $H\in\H^2(\R^d)$, then $\int_0^tU_sH_s\cdot \mathrm{d} X_s$ is a martingale.
  Indeed, by applying the Burkholder-Davis-Gundy inequality, we obtain
  \vspace{-2.5mm}
$$
    \E\Big[\sup_{0\leq t\leq T} \Big\vert \int_0^tU_sH_s\cdot \mathrm{d} X_s\Big\vert \Big]\leq C \E\Big[ \big(\int_0^T|U_sH_s|^2 \mathrm{d} s \big)^\frac12\Big]
    \leq C \Vert U\Vert _{\S^2}\Vert H\Vert _{\H^2}<\infty.
    \vspace{-2.5mm}
$$
  More generally, we can show in the same way that if $p>1$ and $q>1$ are conjugate exponents, and $U\in\S^p$ and $H\in\H^q(\R^d)$, then $\int_0^tU_sH_s\cdot \mathrm{d} X_s$ is a martingale.
  }
\end{Remark}

\subsection{Deriving the $\L^2$ Sensitivity of DROCS Problems}
We now turn to the study of the $\L^2$ sensitivity. 
We reformulate the constraint via duality:
\vspace{-2mm}
$$
  \overline{\V}_2(r)
  =
  \inf_{\alpha \in \Ac}\inf_{\tau\in\Tc}\sup_{\beta\in\Hb(\R^d)}\inf_{\gm>0}
  \Big\{
  J_\tau(\alpha,\beta)
  -
  \E^{\lm^\alpha+\beta}
  \Big[
  \int_0^\tau \frac1{4\gm}|\beta_t|^2\,\mathrm{d}t
  \Big]
  +
  \frac{r^2}{4\gm}
  \Big\},
$$
where $J_\tau$ is defined by \eqref{eqdef:constant problems alpha DROS}.
Define
\vspace{0mm}
$$
\tl{\V}_2(r)
:=
\inf_{\gm>0}
\big\{
G(\gm)+\frac{r^2}{4\gm}
\big\},
\,\,
\text{where }
G(\gm)
:=
\inf_{\alpha \in \Ac}\inf_{\tau\in\Tc}\sup_{\beta\in\Hb(\R^d)}
\hat{J}_\gm(\tau, \alpha,\beta)
\vspace{-1mm}
$$
and 
$
\hat{J}_\gm(\tau, \alpha,\beta)
:=
J_\tau(\alpha,\beta)
-
\E^{\lm^\alpha+\beta}
[
\int_0^\tau \frac1{4\gm}|\beta_t|^2\,\mathrm{d}t
].
$

\noindent The proof follows the same strategy as in the control case:
\begin{romanenumerate}
  \item represent $G(\gm)$ as the initial value of a quadratic RBSDE, see Proposition \ref{stopL2prop:RBSDE gm + BSDE wellposed+ representation G(gm)};
  \item study the differentiability of $\gm\mapsto G(\gm)$ at $\gm=0$, see Propositions \ref{stopL2prop:borne U gm Vgm Kgm}, \ref{stopL2prop:tau gm cv vers tau 0}, and \ref{diffat0G};
  \item prove the strong duality relation $\overline{\V}_2(r)=\tl{\V}_2(r)$ by combining Lemma \ref{lem:strong duality} with Proposition \ref{prop:diff of Gtau(gm)};
  \item conclude by applying Lemma \ref{lem:concluding lemma}.
\end{romanenumerate}
The main additional difficulty, compared with the control case, comes from the reflection term and the dependence of the optimal stopping time on $\gm$.
Consider the following generators, for $\alpha\in\Ac$, $\beta\in\Hb(\R^d)$, and $\gm>0$:
\vspace{-2.5mm}
\begin{equation*}
\begin{split}
  \hat{F}^\gm_t(y,z)&:=\sup_{b\in\R^d}\{f_t(y,z)-b\cdot z-\frac1{4\gm}|b|^2\} 
  =f_t(y,z)+\gm |z|^2, \\
  \text{and} \,\, \hat{h}^{\gm,\alpha,\beta}_t(y,z)&:=l^\alpha_t-k^\alpha_ty -(\lm^\alpha_t+\beta_t)\cdot z-\frac1{4\gm}|\beta_t|^2.
\end{split}
\end{equation*}
\vspace{-1mm}

\begin{Proposition}\label{stopL2prop:RBSDE gm + BSDE wellposed+ representation G(gm)}
Let $\xi$ satisfy Assumption \ref{ass:stop} and $(\xi_t)_{t\in[0,T]}\in\L^\infty$. 
\begin{romanenumerate}
  \item \label{exsite quad RBSDE} For all $\gm\geq0$, the ${\rm RBSDE}(\hat{F}^\gm,\xi)$ admits a unique solution $(Y^\gm,Z^\gm,K^\gm)\in\S^\infty\times\Hb(\R^d)\times\mathbb{A}^2$.
  \item \label{exsitence BSDE stop} For all $\tau\in\Tc$, $\alpha \in \Ac$, $\gm>0$, and $\beta\in\Hb(\R^d)$, the ${\rm BSDE}_\tau(\hat{h}^{\gm,\alpha,\beta},\xi)$ admits a unique solution $(S^{\tau,\gm,\alpha,\beta},\pi^{\tau,\gm,\alpha,\beta})$ such that $(S^{\tau,\gm,\alpha,\beta}_{\wedge\tau},\pi^{\tau,\gm,\alpha,\beta}_{\wedge\tau})\in\S^\infty\times\Hb(\R^d)$.
  \item \label{represnetation formula} Defining
  $
  \tau_t^\gm:=\inf\{s\geq t:Y_s^\gm=\xi_s\}\in\Tc$, 
  $
  \hat{\beta}^\gm:=-2\gm Z^\gm\in\Hb(\R^d)$
  and $
  \hat{\alpha}_t^\gm:=\alpha^*(Y_t^\gm,Z_t^\gm),
  $
  where $\alpha^*$ is defined in Assumption \ref{ass:controlLinf cas controle générales}--\ref{existence argmin}. 
For all $\gm>0$,
  \vspace{-2.5mm}
$$
    \essinf_{\substack{\alpha \in \Ac \\ \tau\in\T_t}}
    \esssup_{\beta\in\Hb(\R^d)} S_t^{\tau,\gm,\alpha,\beta}
    =
    S_t^{\tau_t^\gm,\gm,\hat{\alpha}^\gm,\hat{\beta}^\gm}
    =
    Y_t^\gm
    =
    \esssup_{\beta\in\Hb(\R^d)}
    \essinf_{\substack{\alpha \in \Ac \\ \tau\in\T_t}}
    S_t^{\tau,\gm,\alpha,\beta}.
    \vspace{-1.5mm}
$$
  In particular, for all $\gm>0$, $G(\gm)=Y_0^\gm$.
\end{romanenumerate}
\end{Proposition}

\proof
The quadratic generator $\hat{F}^\gm$ and $\xi$ satisfy the assumptions of Corollary 1 of \citeauthor{Kobylanski2002ReflectedBSDE} \cite{Kobylanski2002ReflectedBSDE}. 
Hence, there exists a unique solution $(Y^\gm,Z^\gm,K^\gm)\in\S^\infty\times\Hb(\R^d)\times\mathbb{A}^2$ of ${\rm RBSDE}(\hat{F}^\gm,\xi)$. 
Moreover, Lemma A-1 from \citeauthor{Hu2022QuadraticMeanFieldReflected} \cite{Hu2022QuadraticMeanFieldReflected} implies that $Z^\gm\in\Hb(\R^d)$. 
This proves Proposition \ref{stopL2prop:RBSDE gm + BSDE wellposed+ representation G(gm)}-\ref{exsite quad RBSDE} .

We now prove \ref{stopL2prop:RBSDE gm + BSDE wellposed+ representation G(gm)}-\ref{exsitence BSDE stop}. For every $\tau\in\Tc$, $\alpha\in\Ac$, $\gm>0$, and $\beta\in\Hb(\R^d)$, the generator $\hat{h}^{\gm,\alpha,\beta}$ is linear in $(y,z)$, with bounded coefficient in $y$ and BMO coefficient in $z$. Therefore, by the well-posedness result for linear BSDEs with BMO coefficients, ${\rm BSDE}_\tau(\hat{h}^{\gm,\alpha,\beta},\xi)$ admits a unique solution $(S^{\tau,\gm,\alpha,\beta},\pi^{\tau,\gm,\alpha,\beta})$ such that $(S^{\tau,\gm,\alpha,\beta}_{\wedge\tau},\pi^{\tau,\gm,\alpha,\beta}_{\wedge\tau})\in\S^\infty\times\Hb(\R^d)$.

We next prove \ref{stopL2prop:RBSDE gm + BSDE wellposed+ representation G(gm)}–\ref{represnetation formula}. 
By the well-posedness of ${\rm BSDE}_{\tau_t^\gm}(\hat{h}^{\gm,\hat\alpha^\gm,\hat{\beta}^\gm},\xi)$, for $ s \in [t,\tau_t^\gm]$
\vspace{-2mm}
$$
\mathrm{d}t\otimes\mathrm{d}\P^0\text{-a.e.},\qquad
(Y_s^\gm,Z_s^\gm)
=
\big(S_s^{\tau_t^\gm,\gm,\hat\alpha^\gm,\hat\beta^\gm},
\pi_s^{\tau_t^\gm,\gm,\hat\alpha^\gm,\hat\beta^\gm}\big)
\in
\S^\infty([0,\tau_t^\gm])\times\Hb(\R^d,[0,\tau_t^\gm]).
$$
By the Skorokhod condition and the continuity of $K^\gm\in\mathbb{A}^2$, we have $\mathrm{d}K^\gm=0$ on $[t,\tau_t^\gm]$. 
Furthermore, for all $\alpha \in \Ac$ and $\beta\in\Hb(\R^d)$,
\vspace{-1mm}
$$
\hat{F}^\gm_t(Y_t^\gm,Z_t^\gm)
=
\hat{h}^{\gm,\hat{\alpha}^\gm,\hat\beta^\gm}_t(Y_t^\gm,Z_t^\gm)
\geq 
\hat{h}^{\gm,\hat{\alpha}^\gm,\beta}_t(Y_t^\gm,Z_t^\gm),
\qquad
\mathrm{d}t \times \mathrm{d}\P^0 \textit{ a.e.}
\vspace{-1mm}$$
By the comparison theorem, for $\alpha \in \Ac$ and $\beta\in\Hb(\R^d)$,
$
S_t^{\tau_t^\gm,\gm,\hat{\alpha}^\gm,\hat\beta^\gm}
\geq 
S_t^{\tau_t^\gm,\gm,\hat{\alpha}^\gm,\beta}.
$
Hence,
\begin{equation}\label{eq:lowebounddpp}
Y_t^\gm
\geq
\esssup_{\beta\in\Hb(\R^d)}
\essinf_{\substack{\alpha \in \Ac \\ \tau\in\T_t}}
S_t^{\tau,\gm,\alpha,\beta}.
\vspace{0mm}
\end{equation}
Fix now $\tau\in\T_t$ and $\alpha \in \Ac$, set
$
\Delta Y:=Y^\gm-S^{\tau,\gm,\alpha,\hat\beta^\gm},
$ and $
\Delta Z:=Z^\gm-\pi^{\tau,\gm,\alpha,\hat\beta^\gm}.
$
Then
\vspace{-1mm}
$$
\mathrm{d}(\Delta Y_s)
=
\Delta Z_s\cdot \mathrm{d}X_s
-
g_s(\Delta Y_s,\Delta Z_s)\mathrm{d}s
+
\mathrm{d}K_s^\gm,
\qquad
\Delta Y_\tau=Y_\tau^\gm-\xi_\tau,
\vspace{-1.5mm}
$$
where
$
g_s(y,z)
:=
f_s^{\hat\alpha^\gm}(Y_s^\gm,Z_s^\gm)
-
f_s^\alpha(Y_s^\gm,Z_s^\gm)
-
k_s^\alpha y
-
(\lm_s^\alpha+\hat\beta_s^\gm)\cdot z.
$
Let
\vspace{-2mm}
$$
\Gamma_s^t
:=
\exp\Big(
\int_t^s k_u^\alpha \mathrm{d}u
-
\frac12\int_t^s |\lm_u^\alpha+\hat\beta_u^\gm|^2\mathrm{d}u
+
\int_t^s (\lm_u^\alpha+\hat\beta_u^\gm)\cdot \mathrm{d}X_u
\Big).
\vspace{-2mm}
$$
By It\^o's formula,
\vspace{-2mm}
$$
\mathrm{d}(\Gamma_s^t\Delta Y_s)
=
\big(\Gamma_s^t\Delta Y_s(\lm_s^\alpha+\hat\beta_s^\gm)-\Gamma_s^t\Delta Z_s\big)\cdot \mathrm{d}X_s
-
\Gamma_s^t
\Big(
-\mathrm{d}K_s^\gm
+
\big(
f_s^{\hat\alpha^\gm}(Y_s^\gm,Z_s^\gm)-f_s^\alpha(Y_s^\gm,Z_s^\gm)
\big)\mathrm{d}s
\Big).
\vspace{-2mm}
$$
Using Theorem 7.2.3 from \citeauthor{Jianfen_BSDE} \cite{Jianfen_BSDE}, Remark \ref{rem:S2*H2=mtg}, and the fact that $\Hb(\R^d)\subset\cap_{p>1}\H^p(\R^d)$, we obtain that
$
\int_t^\cdot
(\Gamma_u^t\Delta Y_u(\lm_u^\alpha+\hat\beta_u^\gm)-\Gamma_u^t\Delta Z_u )\cdot \mathrm{d}X_u
$
is a martingale. 
Moreover,
$
f_s^{\hat\alpha^\gm}(Y_s^\gm,Z_s^\gm)-f_s^\alpha(Y_s^\gm,Z_s^\gm)\leq0.
$
Hence, integrating on $[t,\tau]$,
\vspace{0mm}
$$
\Delta Y_t
=
\E_t\Big[
\Gamma_\tau^t(Y_\tau^\gm-\xi_\tau)
+
\int_t^\tau \Gamma_s^t
\Big(
-\mathrm{d}K_s^\gm
+
\big(
f_s^{\hat\alpha^\gm}(Y_s^\gm,Z_s^\gm)-f_s^\alpha(Y_s^\gm,Z_s^\gm)
\big)\mathrm{d}s
\Big)
\Big]
\leq0
\quad \textit{a.s.}
\vspace{0mm}
$$
Therefore, for all $\tau\in\T_t$ and $\alpha\in\Ac$,
$
Y_t^\gm
=
S_t^{\tau_t^\gm,\gm,\hat\alpha^\gm,\hat\beta^\gm}
\leq
S_t^{\tau,\gm,\alpha,\hat\beta^\gm}
\quad \textit{a.s.}
$
Hence,
$
Y_t^\gm
\leq
\essinf_{\substack{\alpha \in \Ac \\ \tau\in\T_t}}
\esssup_{\beta\in\Hb(\R^d)}
S_t^{\tau,\gm,\alpha,\beta}
$ which along with \eqref{eq:lowebounddpp}, yields the desired result.
\ep

\noindent We now establish the differentiability of $Y_0^\gm$ at the origin.
In contrast with the control case, the proof of differentiability requires additional estimates on the reflection process and on the optimal stopping times. 
We first state the differentiability result, and defer the required estimates on $(Y^\gm,Z^\gm,K^\gm)$ as well as the convergence of the stopping times to the end of this subsection.

\begin{Proposition}\label{diffat0G}
Let $\xi$ satisfy Assumption \ref{ass:stop} and $(\xi_t)_{t\in[0,T]}\in\L^\infty$. 
Let $(Y^\gm,Z^\gm,K^\gm)$ be the unique solution of ${\rm RBSDE}(\hat{F}^\gm,\xi)$ defined in Proposition \ref{stopL2prop:RBSDE gm + BSDE wellposed+ representation G(gm)}. 
Let $\tl\tau$ be defined by \eqref{eqdef:tau_opt_mixed}. 
Let $(U,V)$ be the solution of ${\rm BSDE}_{\tl\tau}(g,0)$, where
\vspace{0mm}
$$
g_t(u,v):=k_t^*u+\lm_t^*\cdot v+|Z_t|^2=: \Delta_t^0(u,v)+|Z_t|^2,
\vspace{0mm}
$$
and where $\lm^*$ and $k^*$ are defined in \eqref{eqdef:tau_opt_mixed}. Then
\vspace{0mm}
$$
\frac{Y_0^\gm-Y_0}{\gm}\xrightarrow[\gm\to0]{}U_0.
$$
\end{Proposition}

\proof
Since $Z\in\Hb(\R^d)\subset\H^4(\R^d)$ by \eqref{ineq:HBMOLp}, the linear BSDE ${\rm BSDE}_{\tl\tau}(g,0)$ defining $(U,V)$ is well posed and, by Proposition 4.1.2 of \citeauthor{Jianfen_BSDE} \cite{Jianfen_BSDE}, admits a unique solution in $\S^2\times\H^2(\R^d)$. Define
$
\dU^\gm:=\delta Y^\gm-U=\frac{Y^\gm-Y}{\gm}-U
$ and
$
\dV^\gm:=\delta Z^\gm-V=\frac{Z^\gm-Z}{\gm}-V.
$
For $t\leq\tl\tau$, $K_t=0$ hence 
\vspace{-3.5mm}
$$
\mathrm{d}\dU_t^\gm
=
\dV_t^\gm\cdot \mathrm{d}X_t
-
\Big(
\Delta_t^\gm(\delta Y_t^\gm,\delta Z_t^\gm)
-
\Delta_t^0(U_t,V_t)
+
|Z_t^\gm|^2-|Z_t|^2
\Big)\mathrm{d}t
+
\frac1\gm\,\mathrm{d}K_t^\gm,
\vspace{-3.5mm}
$$
where, for $\gm>0$,
$
\Delta_t^\gm(y,z):=\frac{f_t(Y_t+\gm y,Z_t+\gm z)-f_t(Y_t,Z_t)}{\gm}.
$
Applying It\^o's formula and taking the expectation (the local martingale is a true martingale by Remark \ref{rem:S2*H2=mtg}), we obtain
\vspace{-3mm}
\begin{align*}
\E\Big[
(\dU_t^\gm)^2
+
\int_{t\wedge\tl\tau}^{\tl\tau}|\dV_s^\gm|^2\,\mathrm{d}s
\Big]
=
&\!-2\E\Big[
\int_{t\wedge\tl\tau}^{\tl\tau}
\dU_s^\gm
\Big(
\Delta_s^\gm(\delta Y_s^\gm,\delta Z_s^\gm)
-
\Delta_s^0(U_s,V_s)
+
|Z_s^\gm|^2-|Z_s|^2
\Big)\mathrm{d}s
\Big]
\\
&-
\frac2\gm
\E\Big[
\int_{t\wedge\tl\tau}^{\tl\tau}\dU_s^\gm\,\mathrm{d}K_s^\gm
\Big].
\end{align*}
\vspace{-5mm}

\noindent Let $\tau^\gm$ be defined by \eqref{eqdef:stoppingtime}. Using the Skorokhod condition, the continuity of $K^\gm$, and integration by parts, as in \eqref{eqdef:estimDKint}, we get
\vspace{-2mm}
$$
-\int_{t\wedge\tl\tau}^{\tl\tau}\dU_s^\gm\,\mathrm{d}K_s^\gm
\leq
-\int_{\tau^\gm}^{\tl\tau}K_s^\gm V_s\cdot \mathrm{d}X_s
+
\int_{\tau^\gm}^{\tl\tau}K_s^\gm |g_s(U_s,V_s)|\,\mathrm{d}s.
\vspace{-3mm}
$$
Therefore,
\vspace{-3mm}
\begin{align*}
-\frac2\gm\E\Big[\int_0^{\tl\tau}\dU_t^\gm\,\mathrm{d}K_t^\gm\Big]
&\leq
2\E\Big[\int_{\tau^\gm}^{\tl\tau}\frac{K_t^\gm}{\gm}|g_t(U_t,V_t)|\,\mathrm{d}t\Big]
\\
&\leq
2\E\Big[
\sup_{\tau^\gm\leq t\leq\tl\tau}\frac{K_t^\gm}{\gm}
\int_{\tau^\gm}^{\tl\tau}|g_t(U_t,V_t)|\,\mathrm{d}t
\Big]
\\
&=
2\E\Big[
\sup_{\tau^\gm\leq t\leq\tl\tau}\frac{|\Delta K_t^\gm|}{\gm}
\int_{\tau^\gm}^{\tl\tau}|g_t(U_t,V_t)|\,\mathrm{d}t
\Big]
\\
&\leq
2\E\Big[\frac1{\gm^2}\sup_{0\leq t\leq T}|\Delta K_t^\gm|^2\Big]^{1/2}
\E\Big[\Big(\int_{\tau^\gm}^{\tl\tau}|g_t(U_t,V_t)|\,\mathrm{d}t\Big)^2\Big]^{1/2}.
\end{align*}
\noindent 
Since $Z\in\H^4(\R^d)$ and 
\textit{a.s.},
$\tau^\gm\rightarrow \tl\tau$, Proposition \ref{stopL2prop:tau gm cv vers tau 0} yields
$
\E
[
(\int_{\tau^\gm}^{\tl\tau}|g_t(U_t,V_t)|\,\mathrm{d}t )^2 ]\rightarrow 0.
$
Moreover, by Proposition \ref{stopL2prop:borne U gm Vgm Kgm},
$
\E [\frac1{\gm^2}\sup_{0\leq t\leq T}|\Delta K_t^\gm|^2 ]^{1/2}
$
is bounded. Hence,
\begin{equation}\label{stopL2eq:integrale deltaUgmKgm tend vers 0}
\frac2\gm\E\Big[\int_0^{\tl\tau}\dU_t^\gm\,\mathrm{d}K_t^\gm\Big]\xrightarrow[\gm\to0]{}0.
\end{equation}

\noindent Next, by triangle inequality, 
\vspace{-2mm}
$$
\E\Big[\int_0^{\tl\tau}|\dU_t^\gm|\,\big||Z_t^\gm|^2-|Z_t|^2\big|\,\mathrm{d}t\Big]
\leq
\Big(
\Vert U^\gm\Vert_{\S^2}
+
\Vert U\mathbf{1}_{[0,\tl\tau]}\Vert_{\S^2}
\Big)
\E\Big[
\Big(
\int_0^{\tl\tau}\big||Z_t^\gm|^2-|Z_t|^2\big|\,\mathrm{d}t
\Big)^2
\Big]^{1/2}.
\vspace{-2mm}
$$

\noindent Since $U\in\S^2([0,\tl\tau])$ and Proposition \ref{stopL2prop:borne U gm Vgm Kgm} holds, this implies
\vspace{-2mm}
\begin{equation}\label{stopL2eq:integrale de deltaUgm (Zgm - Z0) carré cv vers 0}
\E\Big[\int_0^{\tl\tau}|\dU_t^\gm|\,\big||Z_t^\gm|^2-|Z_t|^2\big|\,\mathrm{d}t\Big]\xrightarrow[\gm\to0]{}0.
\vspace{-2mm}
\end{equation}
Furthermore, we can show that
$
\E [
\int_{t\wedge\tl\tau}^{\tl\tau}
|\dU_s^\gm|
\,
|
\Delta_s^\gm(\delta Y_s^\gm,\delta Z_s^\gm)-\Delta_s^0(U_s,V_s)
|
\,\mathrm{d}s
]
\rightarrow 0,
$
by similar consideration as in the proof of \ref{thm:Vup=Vdown mixed inf} of Theorem \ref{thm:sensimixedproblem}.
Combining this convergence with \eqref{stopL2eq:integrale de deltaUgm (Zgm - Z0) carré cv vers 0} and \eqref{stopL2eq:integrale deltaUgmKgm tend vers 0}, we obtain
$
(\dU_0^\gm)^2+\E[\int_0^{\tl\tau}|\dV_t^\gm|^2\,\mathrm{d}t]\rightarrow 0.
$
Therefore,
$
\frac{Y_0^\gm-Y_0}{\gm}\xrightarrow[\gm\to0]{}U_0.
$
\ep

\begin{Remark}\label{stopL2rem:ccl partie asymptotique de G(gm)}
{\rm
By Proposition \ref{stopL2prop:RBSDE gm + BSDE wellposed+ representation G(gm)}, for $\gm>0$, $G(\gm)$ is represented by the solution of an RBSDE at time $0$, namely
$
G(\gm)=Y_0^\gm.
$
Moreover, $G$ can be continuously extended at $\gm=0$ by setting $G(0)=Y_0$, and $G$ is nondecreasing, as follows from the comparison argument used in the proof of Proposition \ref{stopL2prop:tau gm cv vers tau 0}. 
We also proved in Proposition \ref{diffat0G} that $G$ is differentiable at $0$, with
$
G'(0)=U_0.
$
}
\end{Remark}
\vspace{-1mm}
\noindent We now show that $\overline{\V}_2(r)=\tl \V_2(r)$. For this subsection, fix $\tau\in\Tc$, $\alpha\in\Ac$, and $r>0$. 
Consider the following Lagrangian:
\vspace{-4mm}
$$
  \Lc_r^{\alpha,\tau}(\gm,\beta)
  :=
  \E^{\lm^\alpha+\beta}
  \Big[
  \Kc_\tau^\alpha \xi_\tau
  +
  \int_0^\tau
  \Big(
  \Kc_t^\alpha l_t^\alpha
  -
  \frac1{4\gm}|\beta_t|^2
  \Big)\mathrm{d}t
  \Big]
  +
  \frac{r^2}{4\gm},
\vspace{-4mm}
$$
and define
\vspace{-4mm}
$$
  G^{\alpha,\tau}(\gm)
  :=
  \sup_{\beta\in\Hbtau(\R^d)}
  \E^{\lm^\alpha+\beta}
  \Big[
  \Kc_\tau^\alpha \xi_\tau
  +
  \int_0^\tau
  \Big(
  \Kc_t^\alpha l_t^\alpha
  -
  \frac1{4\gm}|\beta_t|^2
  \Big)\mathrm{d}t
  \Big].
  \vspace{-3mm}
$$

\begin{Proposition}\label{prop:diff of Gtau(gm)}
Let $\xi$ satisfy Assumption \ref{ass:stop} and $(\xi_t)_{t\in[0,T]}\in\L^\infty$. 
Then $G^{\tau,\alpha}$ is a nondecreasing continuous function, which admits a continuous extension at $0$. Furthermore, $G^{\tau,\alpha}$ is differentiable on $\R_+^*$, with
$$
(G^{\tau,\alpha})'(\gm)=U_0^{\tau,\alpha,\gm},
$$
where $(U^{\tau,\alpha,\gm},V^{\tau,\alpha,\gm})$ is the solution of the linear ${\rm BSDE}_\tau(g,0)$, with
$$
g_t(u,v)
=
k_t^\alpha u
+
\lm_t^\alpha\cdot v
+
|\pi_t^{\tau,\alpha,\gm}|^2
-
2\gm\,\pi_t^{\tau,\alpha,\gm}\cdot v.
$$
\end{Proposition}

\proof
Following the same line of argument as in the proof of Proposition \ref{prop:diff of Galpha(gm)}, we obtain
$
 G^{\tau,\alpha}(\gm)=S_0^{\tau,\alpha,\gm},
$
where $(S^{\tau,\alpha,\gm},\pi^{\tau,\alpha,\gm})$ is the unique solution of ${\rm BSDE}_\tau(\hat{F}^{\alpha,\gm},\xi)$, with
$
\hat{F}_t^{\alpha,\gm}(y,z)
=
l_t^\alpha-k_t^\alpha y-\lm_t^\alpha\cdot z+\gm|z|^2.
$
This BSDE is well posed by Theorem 6.3.1 from \citeauthor{Jianfen_BSDE} \cite{Jianfen_BSDE}. 
Let $(S^{\tau,\alpha,\gm},\pi^{\tau,\alpha,\gm})$ denote its unique solution, and define
$
\hat{\beta}^{\tau,\alpha,\gm}:=-2\gm\,\pi^{\tau,\alpha,\gm}\in\Hbtau(\R^d).
$
The remainder of the proof follows exactly the same line of argument as in the proof of Proposition \ref{prop:diff of Galpha(gm)}.
\ep

\noindent We can now conclude exactly as in the proof of Theorem \ref{thm:V'(0) Linf et L2}, \ref{thm:V'(0) L2}. Applying Lemmas \ref{lem:strong duality} and \ref{lem:concluding lemma} yields the desired result.

\subsubsection*{Estimates on $(Y^\gm, Z^\gm, K^\gm)$ and convergence of stopping times.}

\begin{Proposition}\label{stopL2prop:borne U gm Vgm Kgm}
Let $\xi$ satisfy Assumption \ref{ass:stop} and $(\xi_t)_{t\in[0,T]}\in\L^\infty$. For $\gm>0$, define
$
U^\gm:=\frac{Y^\gm-Y}{\gm}\in\S^\infty
$, $
V^\gm:=\frac{Z^\gm-Z}{\gm}\in\Hb(\R^d)
$ and $
\Delta K^\gm:=K^\gm-K\in\S^2
$
.
There exists $M>0$ such that, for all $0<\gm\leq1$,
$$
\E\Big[\sup_{0\leq t\leq T}|U_t^\gm|^2+\int_0^T|V_t^\gm|^2\,\mathrm{d}t+\frac1{\gm^2}\sup_{0\leq t\leq T}|\Delta K_t^\gm|^2\Big]\leq M.
$$
\end{Proposition}
\proof
{\bf Step 1.}
We first show that, for every $\gm>0$,
$
\E [\int_0^T|V_t^\gm|^2\,\mathrm{d}t ]\leq \Vert Z^\gm\Vert_{\H^4}^4.
$
First notice that $(U^\gm,V^\gm,\Delta K^\gm)$ satisfies
\begin{equation}\label{stopL2eq:dU gm}
\mathrm{d}U_t^\gm=V_t^\gm\cdot \mathrm{d}X_t-\big(g_t^\gm(U_t^\gm,V_t^\gm)+|Z_t^\gm|^2\big)\mathrm{d}t+\frac1\gm\,\mathrm{d}\Delta K_t^\gm, \,\, U_T^\gm=0,
\end{equation}
where
$
g_t^\gm(u,v):=\frac{f_t(Y_t+\gm u,Z_t+\gm v)-f_t(Y_t,Z_t)}{\gm}.
$
Applying It\^o's formula to $(U^\gm)^2$ yields
$$
\mathrm{d}(U_t^\gm)^2
=
2U_t^\gm V_t^\gm\cdot \mathrm{d}X_t
-
2U_t^\gm\big(g_t^\gm(U_t^\gm,V_t^\gm)+|Z_t^\gm|^2\big)\mathrm{d}t
+
|V_t^\gm|^2\mathrm{d}t
+
\frac2\gm U_t^\gm\,\mathrm{d}\Delta K_t^\gm.
$$
Integrating, we get
\begin{equation}\label{stopL2eq:Ugm carré}
\begin{split}
(U_t^\gm)^2+\int_t^T|V_s^\gm|^2\,\mathrm{d}s
=
-2\int_t^T U_s^\gm V_s^\gm\cdot \mathrm{d}X_s
&+
2\int_t^T U_s^\gm\big(g_s^\gm(U_s^\gm,V_s^\gm)+|Z_s^\gm|^2\big)\mathrm{d}s
\\
&-
\frac2\gm\int_t^T U_s^\gm\,\mathrm{d}\Delta K_s^\gm.
\end{split}
\end{equation}

\noindent Notice that
$
-\int_t^T U_s^\gm\,\mathrm{d}\Delta K_s^\gm\leq0.
$
Indeed, using the Skorokhod condition, we have
\begin{equation*}
\begin{split}
-\int_t^T U_s^\gm\,\mathrm{d}\Delta K_s^\gm
&=
\frac1\gm\int_t^T (Y_s-Y_s^\gm)\,\mathrm{d}K_s^\gm
-
\frac1\gm\int_t^T (Y_s-Y_s^\gm)\,\mathrm{d}K_s
\\
&=
\frac1\gm\int_t^T \underbrace{(Y_s-\xi_s)}_{\leq0}\underbrace{\mathrm{d}K_s^\gm}_{\geq0}
-
\frac1\gm\int_t^T \underbrace{(\xi_s-Y_s^\gm)}_{\geq0}\underbrace{\mathrm{d}K_s}_{\geq0}
\leq0.
\end{split}
\end{equation*}

\noindent 
Therefore, taking expectations in \eqref{stopL2eq:Ugm carré}, and using the Lipschitz assumption on $f$ from Assumption \ref{ass:controlL2 générales}, we obtain
\vspace{-3mm}
\begin{align*}
\E\Big[(U_t^\gm)^2+\int_t^T|V_s^\gm|^2\,\mathrm{d}s\Big]
&\leq
2\E\Big[\int_t^T |U_s^\gm|\big(g_s^\gm(U_s^\gm,V_s^\gm)+|Z_s^\gm|^2\big)\mathrm{d}s\Big]
\\
&\leq
2C\E\Big[\int_t^T |U_s^\gm|(|U_s^\gm|+|V_s^\gm|)\,\mathrm{d}s\Big]
+
2\E\Big[\int_t^T |U_s^\gm||Z_s^\gm|^2\,\mathrm{d}s\Big]
\\
&\leq
\! \Big(2C \!+\! \frac2\eps\Big)\E\Big[\int_t^T \!\! |U_s^\gm|^2\,\mathrm{d}s\Big]
\!+\!
\eps\E\Big[\int_t^T \!\!|V_s^\gm|^2\,\mathrm{d}s\Big]
\!+\!
\eps\E\Big[\int_t^T \!\!|Z_s^\gm|^4\,\mathrm{d}s\Big].
\end{align*}
\vspace{-5mm}

\noindent 
Choosing $\eps=1/2$ and applying the Gronwall inequality, there exists $\tilde C>0$ such that
\vspace{-2mm}
\begin{equation}\label{eqdef:estim1}
\E\Big[(U_t^\gm)^2+\int_t^T|V_s^\gm|^2\,\mathrm{d}s\Big]
\leq
\tilde C\,\E\Big[\int_t^T|Z_s^\gm|^4\,\mathrm{d}s\Big].
\vspace{-2mm}
\end{equation}

\noindent {\bf Step 2.}
We next show that
\vspace{-5mm}
\begin{equation}\label{ineq:step2}
\E\Big[\sup_{0\leq t\leq T}|U_t^\gm|^2\Big]\leq C\Vert Z^\gm\Vert_{\H^4}^4.
\vspace{-3mm}
\end{equation}
Taking the supremum then the expectation in \eqref{stopL2eq:Ugm carré}, by BDG inequality, we get
\vspace{-2.5mm}
\begin{align*}
\E\Big[\sup_{0\leq t\leq T}|U_t^\gm|^2\Big]
&\leq
2\E\Big[\sup_{0\leq t\leq T}\Big|\int_t^T U_s^\gm V_s^\gm\cdot \mathrm{d}X_s\Big|\Big]
+
2\E\Big[\int_0^T |U_s^\gm|\big(g_s^\gm(U_s^\gm,V_s^\gm)+|Z_s^\gm|^2\big)\mathrm{d}s\Big]
\\
&\leq
C\E\Big[\Big(\int_0^T |U_s^\gm|^2|V_s^\gm|^2\,\mathrm{d}s\Big)^{1/2}\Big]
+
2\E\Big[\sup_{0\leq t\leq T}|U_t^\gm|\int_0^T |Z_s^\gm|^2\,\mathrm{d}s\Big]
\\
&\quad+
2\E\Big[\int_0^T |U_s^\gm|(|U_s^\gm|+|V_s^\gm|)\,\mathrm{d}s\Big].
\end{align*}
\vspace{-5mm}

\noindent Using \eqref{eqdef:estim1}, up to a change of $\tilde C$,
\vspace{-3mm}
\begin{align*}
\E\Big[\sup_{0\leq t\leq T}|U_t^\gm|^2\Big]
&\leq
\tilde C\,\E\Big[\int_0^T |Z_s^\gm|^4\,\mathrm{d}s\Big]
+
C\E\Big[\sup_{0\leq t\leq T}|U_t^\gm|\Big(\int_0^T |V_s^\gm|^2\,\mathrm{d}s\Big)^{1/2}\Big]
\\
&\quad+
2\E\Big[\sup_{0\leq t\leq T}|U_t^\gm|\int_0^T |Z_s^\gm|^2\,\mathrm{d}s\Big]
\\
&\leq
\frac14\E\Big[\sup_{0\leq t\leq T}|U_t^\gm|^2\Big]
+
C^2\E\Big[\int_0^T |V_s^\gm|^2\,\mathrm{d}s\Big]
+
\frac14\E\Big[\sup_{0\leq t\leq T}|U_t^\gm|^2\Big]
+
\tilde C\,\Vert Z^\gm\Vert_{\H^4}^4.
\end{align*}
\vspace{-5mm}

\noindent Using again \eqref{eqdef:estim1},
\vspace{-3mm}
\begin{align*}
\E\Big[\sup_{0\leq t\leq T}|U_t^\gm|^2\Big]
&\leq
2C^2\E\Big[\int_0^T |V_s^\gm|^2\,\mathrm{d}s\Big]
+
8\Vert Z^\gm\Vert_{\H^4}^4
\\
&\leq
2C^2\eps\E\Big[\sup_{0\leq t\leq T}|U_t^\gm|^2\Big]
+
\Big(\frac{2C^2}{\eps}+8\Big)\Vert Z^\gm\Vert_{\H^4}^4,
\end{align*}
\vspace{-5mm}

\noindent 
which, combined with \eqref{eqdef:estim1}, yields \eqref{ineq:step2} after choosing $\eps>0$ sufficiently small.

\noindent {\bf Step 3.}
We prove that there exists $M>0$ such that
\vspace{-2mm}
\begin{equation}\label{ineq:step3}
\text{for }0\leq\gm\leq1,\qquad \Vert Z^\gm\Vert_{\Hb}\leq M.
\vspace{-2mm}
\end{equation}
Using Theorem 1 from \citeauthor{Kobylanski2002ReflectedBSDE} \cite{Kobylanski2002ReflectedBSDE}, we get $\Vert Y^\gm\Vert_{\S^\infty}\leq A'$, with $A'$ independent of $\gm\leq1$. 
Furthermore, by the estimate given in the proof of Lemma A-1 from \citeauthor{Hu2022QuadraticMeanFieldReflected} \cite{Hu2022QuadraticMeanFieldReflected}, we obtain
$
\Vert Z^\gm\Vert_{\Hb}
\leq
\frac14\exp(8\Vert Y^\gm\Vert_{\S^\infty})
\leq
\frac{e^{8A'}}4,
$
which completes this step.

\smallskip

\noindent {\bf Step 4.} We conclude.
By Step 3 and \eqref{ineq:HBMOLp}, there exists $K>0$ such that, for all $0\leq\gm\leq1$,
\vspace{-2mm}
\begin{equation}\label{stopL2eq:majo Z gm H4}
\Vert Z^\gm\Vert_{\H^4}^4\leq K.
\vspace{-2mm}
\end{equation}
Combining \eqref{stopL2eq:majo Z gm H4}, \eqref{ineq:step2}, and \eqref{ineq:step3}, we get, for some $M>0$ independent of $\gm\in(0,1]$,
\vspace{-2.5mm}
\begin{equation}\label{stopL2eq:majo U gm V gm intermediaire}
\E\Big[\sup_{0\leq t\leq T}|U_t^\gm|^2+\int_0^T|V_t^\gm|^2\,\mathrm{d}t\Big]\leq M.
\vspace{-2.5mm}
\end{equation}
Integrating \eqref{stopL2eq:dU gm} on $[0,t]$ yields
\vspace{-2mm}
$$
\frac1\gm\Delta K_t^\gm
=
-\int_0^t V_s^\gm\cdot \mathrm{d}X_s
+
\int_0^t \big(g_s^\gm(U_s^\gm,V_s^\gm)+|Z_s^\gm|^2\big)\mathrm{d}s
+
U_t^\gm-U_0^\gm.
\vspace{-4mm}
$$
Therefore,
\vspace{-3mm}
\begin{align*}
\frac1{\gm^2}\sup_{0\leq t\leq T}|\Delta K_t^\gm|^2
&\leq
3^2\sup_{0\leq t\leq T}\Big|\int_0^t V_s^\gm\cdot \mathrm{d}X_s\Big|^2
+
3^2\Big(\int_0^T |g_s^\gm(U_s^\gm,V_s^\gm)|+|Z_s^\gm|^2\,\mathrm{d}s\Big)^2
\\
&\quad+
2\cdot 3^2\sup_{0\leq t\leq T}|U_t^\gm|^2.
\end{align*}
\vspace{-5mm}

\noindent Taking expectations and using the BDG inequality, together with \eqref{stopL2eq:majo Z gm H4}, \eqref{stopL2eq:majo U gm V gm intermediaire}, and the Lipschitz property of $f$, yields the claim.
\ep

\noindent A direct consequence of Proposition \ref{stopL2prop:borne U gm Vgm Kgm} is the following.

\begin{Proposition}\label{stopL2prop:tau gm cv vers tau 0}
Let $\xi$ satisfy Assumption \ref{ass:stop} and $(\xi_t)_{t\in[0,T]}\in\L^\infty$. Define
\vspace{-2mm}
\begin{equation}\label{eqdef:stoppingtime}
\tau^\gm:=\inf\{t\geq0:Y_t^\gm=\xi_t\}\in\Tc,
\qquad \text{for }\gm\geq0.
\vspace{-2mm}
\end{equation}
Then $(\tau^\gm)_{\gm\geq0}$ is decreasing \textit{a.s.} and $\tau^\gm\rightarrow\tl\tau$ \textit{a.s.}
\end{Proposition}

\proof
By the comparison theorem for quadratic RBSDEs from Proposition 3.2 of \citeauthor{Kobylanski2002ReflectedBSDE} \cite{Kobylanski2002ReflectedBSDE}, and since for $0\leq\gm\leq\gm'$,
$
\hat{F}_t^\gm(y,z)\leq \hat{F}_t^{\gm'}(y,z),
$
we have
$
Y_t^\gm\leq Y_t^{\gm'}
\textit{a.s.}
$
Also, by Proposition \ref{stopL2prop:borne U gm Vgm Kgm} 
$
\E [\sup_{0\leq t\leq T}|Y_t^\gm-Y_t|^2 ]
\leq
M\gm^2
\rightarrow 0.
$
The remainder of the argument proceeds as in Step 2 of the proof of \ref{thm:Vup=Vdown mixed inf} of Theorem \ref{thm:sensimixedproblem}.
\ep

\section{Appendix}\label{subsect:DeterministicLemma}

\textit{Proof of Lemma \ref{lem:concluding lemma}.}
By assumption,
$
g(\gm)=g(0)+g'(0)\gm+o(\gm).
$
First, notice that, since $g$ is nondecreasing and continuous at $\gm=0$, we have $V(0)=g(0)$.
For $r>0$, define
\vspace{-4mm}
$$
U(r)
:=
\frac{V(r)-V(0)}{r}
=
\inf_{\gm>0}\Big\{\frac{g(\gm)-g(0)}{r}+\frac{r}{\gm}\Big\}
=
\inf_{\gm>0}\Big\{\frac{g(r\gm)-g(0)}{r}+\frac{1}{\gm}\Big\}\geq0,
\vspace{-2mm}
$$
where the last inequality follows from the monotonicity of $g$. In the following we denote by $\overline{\lim} $ (resp $\underline{\lim}$) the $\limsup$ (resp $\liminf$).
Using the classical inequality ``$(\overline{\lim})\inf\leq\inf(\overline{\lim})$'', we get
\vspace{-5mm}
$$
\mathop{\overline{\lim}}\limits_{r\to 0}
U(r)
\leq
\inf_{\gm>0}
\mathop{\overline{\lim}}\limits_{r\to 0}
 \Big\{\frac{g(r\gm)-g(0)}{r}+\frac{1}{\gm} \Big\}
=
\inf_{\gm>0}\Big\{l\gm+\frac{1}{\gm}\Big\},
\vspace{-2.5mm}
$$
where $l:=g'(0)$.
The problem
$
\inf_{\gm>0} \{l\gm+\frac{1}{\gm} \}
$
has a unique optimizer at $\gm=\frac{1}{\sqrt{l}}$.
Hence,
\vspace{-3mm}
\begin{equation*}
\label{eq:limsup U(r) leq 2sqrt(l)}
\mathop{\overline{\lim}}\limits_{r\to 0}U(r)\leq2\sqrt{l}.
\vspace{-3mm}
\end{equation*}
Since $U\geq0$, if $l=0$, we immediately have $U(r)\to0$, and therefore $V'(0)=0=2 \sqrt{l}$.

\noindent Suppose now that $l>0$.
Consider $(\gm_r)_{r>0}$, a family of $r$-optimizers for $U(r)$, \textit{i.e.},
\vspace{-3mm}
\begin{align}
\label{eq:famille de r optimiseurs}
\frac{g(r\gm_r)-g(0)}{r}+\frac{1}{\gm_r}-r
\leq
U(r)
\leq
\frac{g(r\gm_r)-g(0)}{r}+\frac{1}{\gm_r}.
\end{align}
\vspace{-6mm}

\noindent 
Let $(r_n)_n$ be a sequence of positive numbers converging to $0$.
Since $g$ is nondecreasing, Inequality \eqref{eq:famille de r optimiseurs} implies that
$
0\leq\frac{1}{\gm_{r_n}}\leq U(r_n)+r_n.
$
Therefore, by \eqref{eq:limsup U(r) leq 2sqrt(l)}, the sequence
$
(\frac{1}{\gm_{r_n}})_n
$
is bounded.
Consider therefore a convergent subsequence, still denoted by
$
(\frac{1}{\gm_{r_n}})_n
$,
and let $\lm\geq0$ be its limit.
Using again \eqref{eq:famille de r optimiseurs}, we get
\vspace{-2mm}
$$
0\leq\frac{g(r_n\gm_{r_n})-g(0)}{r_n}\leq U(r_n)+r_n-\frac{1}{\gm_{r_n}}.
\vspace{-4mm}
$$
Therefore,
\vspace{-2mm}
\begin{equation}
\label{eq:g(rn gmn)-g(0)=O(rn)}
g(r_n\gm_{r_n})-g(0)=O(r_n)=o(1).
\vspace{-2mm}
\end{equation}
Suppose that
$
\underline{\lim} r_n\gm_{r_n}=c>0.
$
Then, for $n$ large enough, $r_n\gm_{r_n}\geq c/2$. Hence, using \eqref{eq:g(rn gmn)-g(0)=O(rn)},
$
g(0)=\lim g(r_n\gm_{r_n})\geq g(c/2).
$
Since $g$ is nondecreasing, this implies that
$
g(\gm)=g(0)
$
for all $\gm\in[0,c/2]$, and therefore $g'(0)=l=0$, which contradicts the assumption $l>0$.
Hence,
$
\liminf r_n\gm_{r_n}=0,
$
so there exists a subsequence $\varphi$ such that
$
r_{\varphi(n)}\gm_{r_{\varphi(n)}}\xrightarrow[n\to\infty]{}0.
$
We may therefore use the asymptotic expansion of $g$:
\vspace{-2mm}
$$
\frac{g(r_{\varphi(n)}\gm_{r_{\varphi(n)}})-g(0)}{r_{\varphi(n)}}
=
\gm_{r_{\varphi(n)}}\,l+o(\gm_{r_{\varphi(n)}}).
\vspace{-2mm}
$$
Since $l>0$, this, together with \eqref{eq:g(rn gmn)-g(0)=O(rn)}, implies that
$
\lim \gm_{r_{\varphi(n)}}\neq\infty,
$
and therefore
$
\lm :=\lim\frac{1}{\gm_{r_n}}>0.
$
Using \eqref{eq:limsup U(r) leq 2sqrt(l)} and \eqref{eq:famille de r optimiseurs}, we obtain
\vspace{-3mm}
$$
\inf_{\gm>0}
\Big\{\frac{l}{\gm}+\gm \Big\}
=
2\sqrt{l}
\geq
\overline{\lim} U(r_{\varphi(n)})
\geq
\overline{\lim} 
\Big\{
\frac{g(r_{\varphi(n)}\gm_{r_{\varphi(n)}})-g(0)}{r_{\varphi(n)}}
+
\frac1{\gm_{r_{\varphi(n)}}}
-
r_{\varphi(n)}
\Big\}
=
\frac{l}{\lm}+\lm.
\vspace{-3mm}
$$
Since
$
\inf_{\gm>0}\left\{\tfrac{l}{\gm}+\gm\right\}
$
has $\sqrt{l}$ as its unique optimizer, we necessarily have
$
\lm=\sqrt{l}.
$
Thus, every convergent subsequence of
$
(\tfrac1{\gm_{r_n}})_n
$
has limit $\sqrt{l}$, and therefore
$
\frac1{\gm_{r_n}}\xrightarrow[n\to\infty]{}\sqrt{l}.
$
Finally, applying \eqref{eq:limsup U(r) leq 2sqrt(l)} and \eqref{eq:famille de r optimiseurs}, we deduce
\vspace{-3mm}
$$
2\sqrt{l}
\geq
\mathop{\overline{\lim}} U(r_n)
\geq
\underline{\lim} U(r_n)
\geq
\underline{\lim}
\Big\{
\frac{g(r_n\gm_{r_n})-g(0)}{r_n}
+
\frac1{\gm_{r_n}}
-
r_n
\Big\}
=
2\sqrt{l}.
\vspace{-4mm}
$$
Therefore, for every sequence $r_n\to0$, we have
$
U(r_n)\xrightarrow[n\to\infty]{}2\sqrt{l}.$
\ep

\smallskip

\noindent \textit{Proof of Lemma \ref{lem:strong duality}.}
Define the Lagrangian
$
\Lc(x,\lm):=\psi(x)-\lm\varphi(x)+\lm r.
$
By weak duality, we have
\vspace{-2mm}
$$
\sup_{x \vert \varphi(x) \leq r} \psi(x)
=
\sup_{x \in \Xc}\inf_{\lm>0}\Lc(x,\lm)
\leq
\inf_{\lm>0}\sup_{x \in \Xc}\Lc(x,\lm)
=
\inf_{\lm>0}U(\lm,r),
\vspace{-1mm}
$$
where
$
U(\lm,r):=H(\lm)+\lm r.
$
Now, since $\varphi$ is nonnegative, $H$ is nonincreasing.
Moreover, since $\psi$ is bounded above by $C$ and $\varphi$ is nonnegative, we have that for all $\lm > 0$, 
$
H(\lm)\leq C.
$
Hence, $H$ can be continuously extended at $\lm=0$.
From now on, we do not distinguish $H$ and its extension.
Since $H$ is bounded and $r>0$, when $\lambda \to \infty$ we have
$
U(\lm,r)\rightarrow +\infty.
$
Therefore, since $H$ is continuous on $\R_+$ and convex, there exists $\lm^*\in\R_+$ such that
$
\inf_{\lm>0}U(\lm,r)=U(\lm^*,r).
$

\noindent {\bf Case $1$, $\lm^*>0$.}
In this case, by differentiability and \eqref{eqde:familyoptim}, we get
$
0=r+H'(\lm^*)=r-\varphi(x_{\lm^*}).
$
Hence
$
\varphi(x_{\lm^*})=r
$ and
since
$
H(\lm^*)=\psi(x_{\lm^*})-\lm^*\varphi(x_{\lm^*}),
$
we obtain
\vspace{-2mm}
$$
\psi(x_{\lm^*})
=
H(\lm^*)+\lm^*r
=
\inf_{\lm>0}\big(H(\lm)+\lm r\big)
\geq
\sup_{x : \, \varphi(x) \leq r} \psi(x)
\geq
\psi(x_{\lm^*}).
\vspace{-3mm}
$$

\noindent {\bf Case $2$, $\lm^*=0$.}
In this case, since $\lm\mapsto U(\lm,r)$ is convex, the fact that its minimum is attained at $0$ implies that $\lm\mapsto U(\lm,r)$ is nondecreasing on $\R_+$.
Hence, for all $\lm>0$,
$
H'(\lm)+r\geq0,
$
which means that
\vspace{-1.5mm}
\begin{equation}\label{eq:xlambda_satisfy_constraint}
r\geq \varphi(x_\lm).
\vspace{-1.5mm}
\end{equation}
Now, since the minimizer is $0$,
$
U(\lm,r)=\psi(x_\lm)-\lm\varphi(x_\lm)+r\lm\geq U(0,r)=H(0),
$
and since $\varphi(x_\lm)\geq0$, we get
$
\psi(x_\lm)\geq H(0)+\lm\varphi(x_\lm)-r\lm\geq H(0)-r\lm.
$
Hence, using \eqref{eq:xlambda_satisfy_constraint}, we have
\vspace{-1.5mm}
$$
H(0)
=
\inf_{\lm>0}U(\lm,r)
\geq
\sup_{ \varphi(x) \leq r} \psi(x)
\geq
\psi(x_\lm)
\geq
H(0)-r\lm.
\vspace{-1.5mm}
$$
Since $\lm$ is arbitrary, letting $\lm\to0$ gives
$
H(0)
=
\inf_{\lm>0}U(\lm,r)
=
\sup_{  \varphi(x) \leq r} \psi(x)
=
H(0),
$
which proves the strong duality.
\ep
\vspace{-3mm}

\section*{Declarations}

\paragraph{Data availability.}

Data sharing is not applicable to this article as no datasets were generated or analyzed during the current study.

\paragraph{Conflict of interest.}

The authors declare that they have no conflict of interest.

\printbibliography[heading=bibintoc, section=0]

\end{document}